\documentclass[final,12pt,twoside]{amsart}
\usepackage{graphicx}
\usepackage{amsmath}
\usepackage{amsthm}
\usepackage{showkeys}
\usepackage{amsmath,amssymb}
\usepackage[titlenumbered,boxed]{algorithm2e}
\title[Generalised sifting in black-box groups]
{Generalised sifting in black-box groups}

\newcounter{algonr}

\newtheorem{theorem}{Theorem}[section]

\newtheorem{lemma}[theorem]{Lemma}
\newtheorem{proposition}[theorem]{Proposition}
\newtheorem{corollary}[theorem]{Corollary}

\theoremstyle{definition}

\newtheorem{definition}[theorem]{Definition}

\author{Sophie Ambrose, Max Neunh\"offer, Cheryl E. Praeger, and Csaba Schneider}\address[Ambrose and Praeger]{School of Mathematics
\& Statistics\\
The University of Western Australia\\ 
35 Stirl\-ing High\-way Crawley\\
Western Australia 6009, Australia}
\address[Neunh\"offer]{Lehrstuhl D f\"ur Mathematik\\
Rheinisch-Westf\"alische Technische Hoch\-schule Aachen\\
Templergraben 64\\
52056 Aachen\\
Germany}
\address[Schneider]{Informatics Laboratory\\
Computer and Automation Research Institute\\
The Hungarian Academy of Sciences\\
1518 Budapest Pf.\ 63\\
Hungary}
\email{sambrose@maths.uwa.edu.au\\max.neunhoeffer@math.rwth-aachen.de
\hfill\break praeger@maths.uwa.edu.au\\csaba.schneider@sztaki.hu\break
WWW:
www.maths.uwa.edu.au/$\sim$sambrose\\www.math.rwth-aachen.de/$\sim$Max.Neunhoeffer\hfill\break www.maths.uwa.edu.au/$\sim$praeger\\ www.sztaki.hu/$\sim$schneider}
\date{draft typeset \today}

\begin{document}
\bibliographystyle{plain}

\thanks{Published in {\em LMS J.
Comput. Math.} {\bf 8} (2005) 217--250.}
\begin{abstract}
We present a generalisation of the sifting procedure introduced originally 
by Sims for computation with finite permutation groups, and now used for 
many computational procedures for groups, such as membership testing and 
finding group orders. Our procedure is a Monte Carlo algorithm, 
and is presented and analysed in the context of black-box groups. 
It is based on a chain of subsets instead of a subgroup chain. 
Two general versions of the procedure are worked out in detail, 
and applications are given
for membership tests for several of the sporadic simple groups. 

Our major objective was
that the procedures could be proved to be Monte Carlo algorithms, 
and their costs computed. In addition we explicitly
determined suitable subset chains for 
six of the sporadic groups, and we implemented the algorithms involving
these chains in the {\sf GAP} computational algebra system.
It turns out that sample implementations perform
well in practice. The implementations will be made available publicly
in the form of a {\sf GAP} package.
\end{abstract}

\maketitle

\catcode`\@=11
\def\newmcodes@{\mathcode`\'"27\mathcode`\*"2A\mathcode`\."613A%
 \mathcode`\-"2D\mathcode`\/"2F\mathcode`\:"603A }
\catcode`\@=\active

\def\proof{\noindent{\it Proof.}\quad}
\def\blackbox{\hfill \vrule height4pt width4pt depth2pt\par \vskip .5cm}
\def\qed{\blackbox}

\def\ex{\vskip .3cm\noindent{\bf Example}\ }

\newcommand{\Om}{\Omega}
\newcommand{\De}{\Delta}
\newcommand{\Ga}{\Gamma}
\newcommand{\al}{\alpha}
\newcommand{\be}{\beta}
\newcommand{\ve}{\varepsilon}
\newcommand{\ga}{\gamma}
\newcommand{\vp}{\varphi}
\newcommand{\lam}{\lambda}
\newcommand{\calP}{{\mathcal{P}}}
\newcommand{\soc}{{\rm soc}}
\newcommand{\Core}{{\rm Core}}
\newcommand{\Sym}{{\rm Sym}}
\newcommand{\Aut}{{\rm Aut}}
\newcommand{\Out}{{\rm Out}}
\newcommand{\Diag}{{\rm Diag}}
\newcommand{\Inn}{{\rm Inn}}
\newcommand{\Hol}{{\rm Hol}}
\newcommand{\Prob}{{\rm Prob}}
\newcommand{\la}{\langle}
\newcommand{\ra}{\rangle}
\newcommand{\Z}{\mathbb Z}

\newcommand{\psl}[2]{\mbox{\sf PSL}_{#1}(#2)}
\newcommand{\pgl}[2]{\mbox{\sf PGL}_{#1}(#2)}
\renewcommand{\sp}[2]{\mbox{\sf Sp}_{#1}(#2)}
\newcommand{\gl}[2]{\mbox{\sf GL}_{#1}(#2)}
\newcommand{\sz}[1]{\mbox{\sf Sz}(#1)}
\newcommand{\psu}[2]{\mbox{\sf PSU}_{#1}(#2)}

\renewcommand{\geq}{\geqslant}
\renewcommand{\leq}{\leqslant}

\newcommand{\calQ}{{\cal Q}}

\newcommand{\calH}{{\cal H}}
\newcommand{\HS}{\ensuremath{\hbox{\sf{HS}}}}
\newcommand{\Ly}{\ensuremath{\hbox{\sf{Ly}}}}

\newcommand{\basicsift}{\mbox{\sc BasicSift}}
\newcommand{\basicsiftrandom}{\mbox{\sc BasicSiftRandom}}
\newcommand{\basicsiftcosetreps}{\mbox{\sc BasicSiftCosetReps}}
\newcommand{\ismember}{\mbox{\sc IsMember}}
\newcommand{\newismember}{\mbox{\sc IsMemberConjugates}}
\newcommand{\randel}{\mbox{\sc RandomElement}}
\newcommand{\fail}{\mbox{\sc Fail}}
\newcommand{\false}{\mbox{\sc False}}
\newcommand{\true}{\mbox{\sc True}}
\newcommand{\prob}[1]{\mathsf{Prob}(#1)}

\section{Introduction}\label{intro}

We generalise a sifting procedure introduced originally by Sims~\cite[Section 4]{Sims70}
(see also \cite[Section 2]{Sims97} and~\cite[Chapter~4]{Seress}) for computation with permutation groups.
Our version is given in the context of black-box groups, and
is based on a chain of subsets rather than a subgroup chain. The essential ingredient
is a scheme for sifting a group element $g$ down a descending  chain
\begin{equation}\label{chain}
G_0=S_0\supset S_1\supset\dots\supset S_k
\end{equation}
of non-empty subsets of a subgroup
$G_0$ of a finite group $G$. The sifting procedure seeks elements
$s_0,\dots,s_{k} \in G_0$ such that,
for each $i<k$, $S_is_i\subseteq S_i$ and $gs_0\dots s_i\in S_{i+1}$; in
addition $gs_0\cdots s_{k-1}s_k=1$, and $s_k$ or its inverse lies in $S_k$.
In many instances the $s_i$ will lie in $S_i$, but this is not required
in general. (Conditions on membership for the $s_i$ are given in
Definition \ref{BasicSiftCond} (c).)

A major objective of this work is to give a careful presentation of a randomised generalised
sifting algorithm with an analysis that proves a guaranteed upper bound on the probability of
failure and provides an estimate of the complexity in terms of the input size. We present our
results in a sequence of steps. This `modular' approach enables us to focus 
in our exposition on the new concepts and methods introduced at each stage.
First we present in Section~\ref{sect:sift} a skeleton version of the generalised sifting 
algorithm {\sc Sift} that involves a sequence of basic modules, namely various versions of a procedure
called {\sc BasicSift},  for which only the input and output
requirements are given explicitly. We
prove in Theorem~\ref{thm:sift} that the algorithm {\sc Sift} is a Las Vegas algorithm.

Next, in Sections~\ref{sect:mc} and~\ref{sect:CL}, we present more details of the versions of {\sc BasicSift}
we have developed, and prove in Theorems~\ref{prop:basic-mc}
and~\ref{prop:basic-cr} that for these versions,
{\sc BasicSift} is a Monte Carlo algorithm. This exposition of {\sc BasicSift} is
given in terms of a generic membership test {\sc IsMember} for which only the input and output
requirements are given explicitly. 
Note that the {\sc BasicSift} modules will often be Monte Carlo algorithms
with a non-zero probability of returning an incorrect result. However
the complete algorithm {\sc Sift} is a Las Vegas algorithm since we can test
with certainty that, for our output element $x=s_0\dots s_k$, the
element $gx$ is equal to the identity. (See Definition \ref{algdefs} for a
discussion of these types of algorithms.)

In Section~\ref{sect:conj}, we introduce a  version of
{\sc IsMember} based on random conjugates. It was this version that inspired 
the development of the conceptual
framework presented in the paper. The idea can best be understood by briefly considering the
following special case. Suppose that a finite group $G$ has a chain of subgroups
\begin{equation}\label{subgroupchain}
G =H_0> H_1>\dots> H_k=\{1_G\}
\end{equation}
and that $a\in H_{k-1}\setminus\{1\}$ is such that, for each $i$, the subset $a^G\cap H_i$ of
$a$-conjugates lying in $H_i$ forms a single $H_i$-conjugacy class $a^{H_i}$. Then for $x\in G$, the
conjugate $a^x$ lies in $H_i$ if and only if $a^x=a^h$ for some $h\in H_i$, and, in turn, this
holds if and only if $xh^{-1}\in C_G(a)$. Thus $a^x\in H_i$ if and only if $x\in C_G(a)H_i$,
that is to say, a membership test for $a^x$ to lie in the subgroup $H_i$ is equivalent to
a membership test for $x$ to lie in the subset $C_G(a)H_i$. Development of this
idea to handle the general case where the subsets $a^G\cap H_i$ split into 
several $H_i$-conjugacy classes led to the theory presented in 
Section~\ref{sect:conj}.

In Section~\ref{sect:orders}, we give full details of
a version of {\sc IsMember} that relies on element orders. For the corresponding version
of {\sc BasicSift} we are then able to provide in Corollary~\ref{final}, our most
comprehensive complexity estimate.

Before presenting the theoretical details we give a worked example of our 
algorithm for the Higman-Sims sporadic simple group in Section~\ref{example}. This 
example  was chosen to illustrate most of the methods that will be developed in 
the paper.

The original motivation for this research stems from the matrix group 
recognition project, see
\cite{LeedhamGreen00,LeedhamGreenNiemeyeretal03}, and
in particular the need to recognize constructively all quasi-simple matrix groups over finite fields.
The usual approach has been to design algorithms for recognizing finite quasi-simple
groups by their intrinsic properties as abstract groups rather than building different algorithms
for each of their different matrix representations. This has resulted
in the development of recognition algorithms for most of the almost simple groups represented
as black-box groups (see \cite{BealsLeedhamGreenetal02,BealsLeedhamGreenetal03,BratusPak00,
CellerLeedhamGreen98,CoopermanFinkelsteinetal97a,KantorSeress01}). A black-box group
is one in which the elements are represented (possibly non-uniquely) as binary strings of
bounded length and in which we can perform the following operations (and only these): we can test
whether two given strings represent the same group element, and
we can produce strings representing the inverse of a given element, and the product of two
given elements. In this paper we give algorithms that involve only these `black-box operations'
of equality tests, extracting inverses, and multiplying group elements. Thus our algorithms are
black-box algorithms.

We are aware of the impressively successful practical algorithms of~\cite{hlo'brw} for recognizing sporadic groups based on the theory of
involution centralisers. However, there seemed to be no framework available 
to analyse the probability of completion or the cost of these algorithms. 
Our motivation was based on both experience and hope: experience with 
developing
recognition algorithms for finite symmetric and alternating groups in 
\cite{BealsLeedhamGreenetal02,BealsLeedhamGreenetal03} complete with proofs
and complexity analyses; and hope that the ideas of Charles Sims could be made 
effective for black-box groups, where information needed about a permutation 
or matrix action must be derived from purely group theoretic properties. 
Success in computing with some of the sporadic simple groups suggested that 
our new approach would provide an alternative method for recognizing and 
computing with these groups. We believe that we have been successful, 
both theoretically and in practice. The algorithmic framework presented
in this paper offers an effective and convenient means of analysing membership 
tests for sporadic simple groups and other groups,
providing proofs of completion probability and complexity. The framework offers
flexibility in choice of subset chains and types of the basic sifting 
procedures. Explicit examples of the algorithms 
have been developed and implemented for several of the sporadic groups and perform very well in 
practice. 
In Section~\ref{sect:appns} we summarise the information about these
examples and also present some details concerning the 
implementations of the procedures presented in this
paper. We emphasise that all groups that occur in this paper are finite.

\section{Generalised sifting: an example}\label{example}

The aim of this section is to explain our approach using the example
of the Higman-Sims group $\HS$. We think of $\HS$ as a group given to
us in its most natural representation, that is, a group of permutations
with degree~100. Throughout this section we use various facts
concerning $\HS$, and the validity of these facts can easily be
checked using the {\sc Atlas}~\cite{Atlas}, or a computer algebra package, such as {\sf
GAP}~\cite{GAP} or {\sc Magma}~\cite{Magma}. In order to describe subgroups
of $\HS$ we use the notation introduced in the {\sc Atlas}.

Suppose that $a, b$ are standard generators in the
sense of \cite{Wilson96} for $\HS$ given on the {\sc Atlas} web
site~\cite{wwwatlas}. Assume that $G$ is a black-box group isomorphic
to $\HS$ and $x, y$ 
are standard generators for $G$
obtained using the procedure described in the online {\sc
Atlas}~\cite{wwwatlas}. Then the map $a\mapsto x,\ b\mapsto y$ can be
extended in a unique way to an isomorphism $\varphi:\HS\rightarrow G$.
Since $\HS$ is a permutation group, it is possible
to compute, using the Schreier-Sims Algorithm, a base and a strong
generating set for $\HS$. Using them, a permutation in $\HS$ can efficiently be
written as a word in $a,\ b$. Thus, if $u\in\HS$ then $\varphi(u)$, as
a word in $x$ and $y$, can be computed efficiently. 
The constructive recognition of the black-box group $G$ requires us to 
perform the opposite process: given $g\in G$, we must find an element
$u\in \HS$ such that $\varphi(u)=g$. This is equivalent to writing 
the element $g$ as a word in $x$ and $y$. 

In order to complete our task, we specify some (precomputed and stored)
elements and subgroups
in $G$. We use the following important convention:
\begin{quotation}
every element we introduce in $G$
from now on will be expressed as a word in $x,\ y$. Similarly, every
subgroup of $G$ we use will be given with a generating set, and each
generator in this set is assumed to be a word in $x,\ y$. 
\end{quotation}

Let $L_1$ be a maximal subgroup of $G$ isomorphic to 
$U_3(5).2$. A generating set for such a subgroup can be found by computing 
a generating set for a maximal subgroup in $\HS$ 
isomorphic to 
$U_3(5).2$, and mapping the
generators into $G$ using $\varphi$. 
In the same way, we find a maximal subgroup $L_2$ in $L_1$
isomorphic to $5^{1+2}:(8:2)$.
Let $L_3$ be a cyclic subgroup of $L_2$ of order $8$ in a complement
$8:2$ for $5^{1+2}$. To be consistent with the notation to be introduced
in later sections of the paper we will denote a generator of $L_3$ by
$a$. We emphasise that this element $a$ lies in $L_3$ and is not a
standard generator of $\HS$.
Let $a$ be
an element of order $8$ in $L_3$, and set $L_4=1$. 
The four generators of $L_3$ are all conjugate to each other in $G$ and in
$L_1$; they fall into two conjugacy classes of $L_2$, and 
they are pairwise not conjugate in $L_3$. 
Thus there are elements $t_1\in L_1$, $t_2,\ t_3,\ t_4\in L_2$ 
such that $a^G\cap L_2=a^{L_2}\cup a^{t_1L_2}$, 
$a^{L_2}\cap L_3=\{a,a^{t_2}\}$ and $a^{t_1L_2}\cap L_3=\{a^{t_1t_3},a^{t_1t_4}\}$. 
Set $\mathcal T_2=\{1,t_1\}$, 
$\mathcal T_3=\{1,t_2,t_1t_3,t_1t_4\}$, and $\mathcal T_4=\{1\}$. 

We therefore have a chain of subgroups $$
G\geq L_1\geq L_2\geq L_3\geq L_4=1,
$$
with $|G:L_1|=176$, $|L_1:L_2|=126$, $|L_2:L_3|=250$, $|L_3|=8$.

\subsection{Sifting $g\in G$ into the first subset: element orders}
Let $g\in G$. If we were to perform Sims's usual sifting procedure, we would
look for an element $h_1\in G$ such that $gh_1\in
L_1$. The probability that a random $h_1$ satisfies this 
property is $|L_1|/|G|=1/176$. What we do instead is as follows. Let
$C_1=C_G(a)$. We look for
an element $h_1\in G$ such that $gh_1\in
C_1L_1$. As $|C_1|=16$ and $|C_1\cap L_1|=|C_{L_1}(a)|=8$, 
the probability that $gh_1\in
C_1L_1$, for a random $h_1$, 
is 
$|C_1L_1|/|G|=2|L_1|/|G|=1/88$.

In order to make this work, we must have a membership
test for $C_1L_1$. Since
$a^G\cap L_1=a^{L_1}$, we have, as explained in the introduction, that, 
for $u\in G$, $u\in C_1L_1$ if and only if $a^u\in L_1$. 
Thus to obtain a membership test for $C_1L_1$, we only need to
design a membership test for $L_1$.
Let $u\in G$, and let $X_1$ be a generating set for $L_1$; set
$\overline{X_1}=X_1\cup\{u\}$. It is clear that $u\in L_1$ if and only if
$\left<\overline{X_1}\right>=\left<X_1\right>$. Now about one
quarter of the elements of $G$ have order~15 or~11, 
but no element in $L_1$ has order equal to one of these numbers.
Hence we select random elements in
$\left<\overline{X_1}\right>$. If such a random element  has order~11
or~15, then we conclude with certainty 
that $u\not\in L_1$. If, however, after many
random selections we do not find an element with order~11 or~15, then
we may say that $u\in L_1$ with a certain 
high probability. This can be formulated to 
give a 
one-sided Monte Carlo membership test for $C_1L_1$; see 
Section~\ref{sect:orders} for details.

\subsection{Sifting $gh_1$ into the second subset: random conjugates}
\label{sect:secondsubset}
The intersection $a^{L_1}\cap L_2$ is the union of two conjugacy
classes in $L_2$, namely $a^{L_2}$ and $a^{t_1L_2}$ where 
$t_1\in L_1$ and we set $\mathcal T_2=\{1,t_1\}$ 
as above. Let $C_2$ denote the set $C_G(a)\mathcal T_2$. As $L_2\leq
L_1$ and $\mathcal T_2\subset L_1$, we have $C_2L_2\subset C_1L_1$. 
Now we seek an element $h_2\in L_1$ such that $gh_1h_2\in
C_2L_2$. 
We will call an element $h_2\in L_1$ `good' if and only if $gh_1h_2\in C_2L_2$, or equivalently, if and only if $a^{gh_1h_2}\in L_2$. 
If $h_2$ is a uniformly distributed random element of $L_1$ and $gh_1\in C_1L_1$, 
then
$a^{gh_1h_2}$ is a uniformly distributed random element of the conjugacy class
$a^{L_1}$. For each $x\in a^{L_1}\cap L_2$ there are $|C_{L_1}(a)|$ choices
for $h_2\in L_1$ such that $a^{gh_1h_2}=x$. Therefore the total number of
`good' elements is $|a^{L_1}\cap L_2||C_{L_1}(a)|$, and so the probability
that $h_2$ is `good' is $|a^{L_1}\cap L_2||C_{L_1}(a)|/|L_1|=|a^{L_1}\cap L_2|/|a^{L_1}|=500/31500=1/63$.

In order to test whether $a^{gh_1h_2}\in L_2$, recall that
$L_2$ is isomorphic to $5^{1+2}:(8:2)$.
A deterministic membership
test for $L_2$ can easily be designed using the fact that
$N_G(Z(5^{1+2}))=L_2$ where $Z(5^{1+2})=\langle b\rangle$ is the centre of 
$5^{1+2}$: namely, to test whether an element $x\in G$ lies in $L_2$
simply test whether $b^x\in \{b,b^2,b^3,b^4\}$.

\subsection{Sifting $gh_1h_2$ into the third subset}
\label{sect:thirdsubset}
The group $L_3$ is cyclic with order~8. 
Set $C_3=C_G(a)\mathcal T_3$ where $\mathcal T_3=\{1,t_2,t_1t_3,t_1t_4\}$. 
As $\mathcal T_3\subset \mathcal T_2L_2$, we obtain $C_3L_3\subset C_2L_2$. We look for an element $h_3\in L_2$ such that,
given $gh_1h_2\in C_2L_2$, we have
$gh_1h_2h_3\in C_3L_3$. Using the definition of $\mathcal T_3$, we
obtain that, given $gh_1h_2\in C_2L_2$,  the condition $gh_1h_2h_3\in C_3L_3$ holds if and only if $a^{gh_1h_2h_3}\in
L_3$. Arguing as for the previous case, the probability that, 
given $gh_1h_2\in C_2L_2$, a random $h_3\in L_2$ yields $gh_1h_2h_3\in C_3L_3$
is at least $\min\left\{|a^{L_2}\cap L_3|/|a^{L_2}|,|a^{t_1L_2}\cap L_3|/|a^{t_1L_2}|\right\}$. It is easy to compute that this number is $2/250=1/125$. 
At the end of this process we have with high probability that
$a^{gh_1h_2h_3}\in\{a,a^{t_2},a^{t_1t_3},a^{t_1t_4}\}$. 
Therefore after a number of equality tests we obtain a word $w$ in $x,\ y$
such that $gw\in C_4$ where $C_4=C_G(a)$. As $|C_G(a)|=16$, using the map $\varphi$, it
is easy to compute each element of $C_G(a)$ as a word in $x,\ y$. Then
comparing $gw$ against the elements of $C_G(a)$, it is now easy to
express $g$ as a word in $x,\ y$.

Thus the main ingredients of this process are a descending chain of
subgroups $\{L_i\}_{i=1}^4$,
a sequence of subsets $\{C_i\}_{i=1}^4$ defined in terms of
the centraliser of the element $a$, and the sequence 
$\{\mathcal T_i\}_{i=1}^4$
of subsets where we take $\mathcal T_1=\{1\}$. Our sifting
procedure progressed through the following descending chain of non-empty
subsets:
$$
G\supset C_1L_1\supset C_2L_2\supset C_3L_3\supset C_4;
$$
the final step was a series of equality tests with the elements of $C_4$.

\section{A small toolbox}\label{probsec}

In this section we collect several results that we need in our proofs.
For an event $E$, $\prob E$ denotes
the probability of $E$. For events $A$ and $B$, $\prob{A | B}$ denotes the probability of $A$, given that $B$ holds. We recall that $\prob{A|B}=\prob{A\cap B}/\prob B$.
The following result from elementary probability
theory will often be used in this article.

\begin{lemma}\label{lem:prob} 
If $A$, $B$, $C$ are events such that $C\subseteq B\subseteq A$,
then $$
\prob{C | A} = \prob{C | B}\cdot\prob{B | A}.
$$
\end{lemma}

\begin{proof}
As $B=B\cap A$ and $C=C\cap B= C\cap A$, we obtain
\begin{multline*}
\prob{C | B}\cdot\prob{B | A}
=\frac{\prob{C\cap B}}{\prob B}\cdot
\frac{\prob{B\cap A}}{\prob A}\\
=\frac{\prob{C\cap A}}{\prob B}\cdot\frac{\prob B}{\prob A}=\frac{\prob{C\cap A}}{\prob A}
=\prob{C | A}.
\end{multline*}
\end{proof}

\begin{lemma}\label{rem:mc}
If $0\leq x<1$, then $\log((1-x)^{-1})\geq x$. 
\end{lemma}
\begin{proof}
Observe that the function
$f(x)=x-\log((1-x)^{-1})$ is strictly decreasing for $0\leq x<1$ and
$f(0)=0$. 
\end{proof}

The following is a general version of Dedekind's modular law. 
Its proof can be carried out following that of~\cite[1.3.14]{rob}.

\begin{lemma}\label{dede}
If $U$ and $V$ are subsets and $Z$ is a subgroup
of a group such that $VZ\subseteq V$ then $(V\cap U)Z=V\cap (UZ)$. 
\end{lemma}

In this paper we use several types of randomised algorithms, that is, algorithms that involve a random choice at some point, so that they do not behave in the same way every time the algorithm is run. 
We also use algorithms which involve no random choices, that is, deterministic algorithms. We collect together here the definitions of these types of algorithms. To aid our exposition we give slightly different definitions of these algorithm types than normal, and we comment on the differences below.

\begin{definition}\label{algdefs}
(a) Let $\ve$ be a real number satisfying $0\leq \ve<1/2$. A {\em Monte Carlo algorithm} with `error probability' $\varepsilon$ is an algorithm that always terminates after a finite number of steps, such that the probability that 
the algorithm gives an incorrect answer is at most $\varepsilon$.

(b)\ A {\em one-sided Monte Carlo} algorithm is a Monte Carlo algorithm which has two types of
output (typically `yes' and `no'), and one of the answers is guaranteed to be correct. 

(c)\ A {\em Las Vegas algorithm} with `failure probability' $\varepsilon$ (where $0\leq\varepsilon<1/2$) terminates after a finite number of steps and either returns an answer, or reports failure. An answer, if given, is always correct, while the probability that the algorithm reports failure is at most $\varepsilon$.

(d) For the purposes of this paper, a {\em deterministic algorithm} is a 
Monte Carlo algorithm for which the `error probability' $\varepsilon$ is $0$, 
or equivalently, a Las Vegas algorithm for which the
`failure probability' $\ve$ is $0$.

\end{definition}

Note that our definitions of Monte Carlo and Las Vegas algorithms vary from the usual ones in that we allow $\varepsilon$ to be zero. The reason for this is that some versions of our \basicsift\ algorithm may be deterministic, that is, have zero probability of failure or of returning an incorrect answer. For ease of exposition we decided to treat such an algorithm as a special case of a Monte Carlo or Las Vegas algorithm.

\section{The generalised sifting algorithm}\label{sect:sift}

In this section we present an algorithm for sifting an element $g$ of
a finite group $G$ down a (given and precomputed) descending chain
(\ref{chain}) of subsets of a subgroup $G_0$ of $G$.
The algorithm returns either {\sc Fail}, or
a word $x=s_0\dots s_{k} \in G_0$ 
such that $gx=1$, $S_is_i\subseteq S_i$
for each $i<k$, and $s_k$ or its inverse lies in $S_k$.
If $g\in G_0$, then (see Theorem~\ref{thm:sift}) the probability that the 
algorithm returns {\sc Fail} is proved to be at most
some pre-assigned quantity $\ve$. Usually the $s_i$ are returned as words 
in a given set $Y$ of generators for $G_0$, or as straight line programs from 
the given generating set $Y$. The algorithm is applied in one of the following 
contexts.
\begin{enumerate}
\item[(1)] The element $g$ is known to lie in $G_0$ and the purpose of the 
algorithm is to express $g$ as a word in a given generating set. 
In this context, 
Theorem~\ref{thm:sift} proves that the algorithm
fails with probability at most $\ve$, for some pre-assigned non-negative real number $\ve<1/2$.
Hence, in this context, Algorithm~\ref{sift} is a \emph{Las 
Vegas algorithm}.
\item[(2)] We only assume that $g\in G$, and the aim is to discover whether or not $g$ lies
in $G_0$. In this context, Theorem~\ref{thm:sift} proves that if 
the algorithm returns an expression for $g$, then $g$ must lie in $G_0$.
On the other hand, if the algorithm returns {\sc Fail} then the element $g$ 
may or may not lie in $G_0$. Moreover, if $g\in G_0$, then 
the probability that the algorithm will return \fail\ is less than some pre-assigned real number $\varepsilon$ where $0\leq\varepsilon< 1/2$.
Hence, in this context (if we interpret the result {\sc Fail} as a finding that $g\not\in G_0$), Algorithm~\ref{sift} is a \emph{one-sided Monte Carlo algorithm}.
\end{enumerate}

In either case we allow the probability bound  $\varepsilon$ to be zero, and in this situation the resulting algorithm
is deterministic.
The basic building block for our algorithm is described in the following definition.

\begin{definition}\label{BasicSiftCond}
A $4$-tuple $(G_0,H,K,\basicsift)$ 
is said to satisfy the {\em basic sift condition} in a group $G$, 
if the following hold:
\begin{enumerate}
\item[(a)] $G$ is a finite group with a subgroup $G_0$;
\item[(b)] $H$ and $K$ are non-empty subsets of $G_0$ such that
either $K=\{1\}$ or $K\subset H$;
\item[(c)] $\basicsift$ is a Monte Carlo algorithm
whose input is a pair $(g,\varepsilon)$, where $g\in G$ and $\varepsilon$
is a non-negative real number. It satisfies the following condition, either
for all inputs $(g,0)$  (in which case it is a deterministic algorithm), or
for all inputs $(g,\varepsilon)$ with $0<\varepsilon<1/2$. \em{The output
$y$ is either {\sc Fail}, or an element of $G_0$ 
such that $Hy\subseteq H$
(if $K\subset H$) or $y^{-1}\in H$ (if $K=\{1\}\not\subset H$).
Moreover, if $g\in H$, then 
$\prob{y=\fail\mbox{, or } (y\in G_0 \mbox{ and }
	gy\not\in K)}\leq \ve.$}
\end{enumerate}
\end{definition}

To avoid confusion we comment on the formulation of the condition in
Definition \ref{BasicSiftCond} (c). Note that $H$ is in general not a
subgroup, and hence $Hy\subseteq H$, for $y\in G$, does not imply that
either of $y$ or $y^{-1}$ lies in $H$. After considering many special cases,
we realised that the set inclusion $Hy\subseteq H$ was the appropriate
requirement.

Suppose that $G$ is a finite group with a subgroup $G_0$ and 
$$
G_0=S_0\supset S_1\supset \cdots\supset S_{k-1}\supset S_k
$$
is a chain of non-empty subsets of $G$, and set $S_{k+1}=\{1\}$. 
Suppose further that, for
$i=0,\ldots,k$, $\basicsift_i$ is an algorithm such that $(G_0,S_{i},S_{i+1},\basicsift_i)$
satisfies the basic sift condition in $G$. 
Then there is a Las Vegas algorithm that, for a given $g\in G$, 
returns either `failure' or an element $s_0s_1\cdots
s_{k}$ of $G_0$ 
such that  $S_is_i\subseteq S_i$ for each $i<k$, the element 
$s_k$ or its inverse lies in $S_k$,
and $gs_0s_1\cdots s_{k}=1$.
Indeed, as shown in Theorem~\ref{thm:sift}, Algorithm~\ref{sift} has
this property.

\begin{algorithm}
\caption{The generalised sift algorithm}
\label{sift}
\SetKw{Input}{Input:}
\SetKw{Output}{Output:}
\SetKwFunction{BasicSift}{\sc BasicSift}
\SetKw{Set}{set}
Algorithm~\ref{sift}: {\sc Sift}\\
\tcc{see Theorem~\ref{thm:sift} for notation}
\Input{\rm $g\in G$ and $(\ve_0,\ldots,\ve_k)$ with $\ve_i\geq0$ and $\sum_i\ve_i<1/2$}\;
\Output{\rm either $x=s_0\cdots s_k$ with $S_is_i\subseteq S_i$
for $i<k$ and $gx=1$, or \fail}\;
\Set $x=1$\;
\For{$i=0$ \KwTo $k$}
	{\Set $s_i=\BasicSift_{i}(gx,\ve_i)$\;
	\eIf{$s_i=\fail$}
	{\Return \fail}
	{\Set $x=xs_i$}}
\eIf{$gx\ne 1$}{\Return \fail}{\Return $x$}
\end{algorithm}

\begin{theorem}\label{thm:sift}
Suppose that $G$, $G_0$, $S_0,\ldots,S_{k+1}$, and $\basicsift_0,\ldots,
\basicsift_k$ are as in the previous paragraph, and let {\sc Sift} denote 
Algorithm~$\ref{sift}$. Let $g\in G$ and $\ve_0,\dots,\ve_k$ be non-negative real numbers such that $\sum_i \ve_i < 1/2$. Then the following hold.
\begin{enumerate}
\item[(i)] If $\mbox{\sc Sift}(g, (\ve_0,\dots,\ve_k))$ 
returns a group element $x$, then $g=x^{-1}\in G_0$ and $x=s_0s_1\cdots s_k$,
where $S_is_i\subseteq S_i$ for each $i\in\{0,1,\dots,k-1\}$, 
and $s_k\in S_k$ if $S_k$ contains $1$, while $s_k^{-1}\in S_k$ otherwise.
\item[(ii)]
The conditional probability that
$\mbox{\sc Sift}(g, (\ve_0,\ldots,\ve_k))$
returns {\sc Fail},  given that $g\in G_0$, is at most $\sum_{i}\ve_i$.
\end{enumerate}
\end{theorem}

\begin{proof}
(i) Suppose that a group element $x=s_0s_1\dots s_k$
is returned. Then the $s_i$ are group elements computed as in 
Algorithm~\ref{sift}. From Definition \ref{BasicSiftCond} (c), since 
each $s_i$ is a group element, we have
that $S_is_i\subseteq S_i$ for each $i\in\{0,1,\dots,k-1\}$, and also for
$i=k$ if $1\in S_k$; while if $1\not\in S_k$, then $s_k^{-1}\in S_k$. Further,
if $1\in S_k$, then $S_ks_k$ contains $s_k$, and hence $S_k$ contains $s_k$.
Finally, for each $i$, $s_i$ lies in $G_0$ since the algorithm $\basicsift_i$
involves random selections from  the group $G_0$. Moreover, by the last
{\bf if} statement of Algorithm~\ref{sift} we have $gx=1$ so that $g=x^{-1}\in G_0$.

(ii) Let $E_0$ denote the event that $g\in G_0$, 
and recall that $G_0=S_0$ and $S_{k+1}=\{1\}$. 
For each 
$i=1,\dots,k$, 
let $E_i$ denote the event that the $i$-th execution of the {\bf for} loop in 
Algorithm~\ref{sift} is attempted, is successful and returns a correct answer. In other words,
\[
E_i:\quad E_{i-1}\mbox{ holds},\ 
S_js_{j}\subseteq S_{j}\ \mbox{for all}\ j=0,\dots,i-1,\ \mbox{and}\ gs_0\dots s_{i-1}\in S_i.
\]
Also define $E_{k+1}$ to be the event that the final execution  of the {\bf for} loop is attempted, is
successful and returns a correct answer. That is,
$$
E_{k+1}:\quad E_k\mbox{ holds},\ S_js_j\subseteq S_j\mbox{ for all }
j=0,\ldots,k-1,\ \mbox{and }gs_0\cdots s_{k}=1.
$$
Then the probability that Algorithm~\ref{sift} returns $x=s_0s_1\dots s_{k}$ with $S_is_i\subseteq S_i$
for all $i=0,\dots,k-1$, and $gx=1$, given that $g\in G_0$, is, by
definition, $\prob{E_{k+1}\,|\,E_0}$.

Now $E_{k+1}\subseteq E_{k}\subseteq \dots\subseteq E_0$, and hence by 
several applications
of Lemma~\ref{lem:prob}, we have that $\prob{E_{k+1}\,|\,E_0}=\prod_{i=0}^{k}
\prob{E_{i+1}\,|\,E_{i}}$. Since $(G_0,S_i,S_{i+1},\basicsift_i)$
satisfies the basic sift condition in $G$ for each $i=0,\dots,k$,
$\prob{E_{i+1}\,|\,E_{i}}\geq 1-\ve_{i}$ for each $i=0,\dots,k$. Hence
\[
\prob{E_{k+1}\,|\,E_0}\geq\prod_{i=0}^{k}(1-\ve_{i}).
\]
Since $0\leq \ve_i<1$ for all $i$, we have 
$\prod_{i}(1-\ve_{i})\geq 1-\sum_{i}\ve_{i}$ (use induction on $k$),
and hence the required probability in part~(ii) is
at most $\sum_{i}\ve_{i}$. \qed
\end{proof}

Algorithm~\ref{sift} allows different types of
algorithms  
to be used for different links of the chain. For example, if 
$S_{k}$ is small, then $\basicsift_k$ relies sometimes on nothing more than an 
exhaustive search through the elements of $S_{k}$ with the parameter $\ve_{k}=0$. Two special types
of $\basicsift$ algorithms are described in detail in Sections~\ref{sect:conj} 
and~\ref{sect:orders}.
We first explore their common
properties as one-sided Monte Carlo algorithms in Sections~\ref{sect:mc} and~\ref{sect:CL}.

\section{{\sc BasicSift}: a general approach}\label{sect:mc}

In this section we present a general approach to designing a $4$-tuple that satisfies
the basic sift condition.
The results of this section will become relevant 
in the discussion of the two algorithms in
Sections~\ref{sect:conj}
and~\ref{sect:orders}. 
We will use one of the general methods given in this section in nearly all
cases when we wish to sift an element of $S_i$ into the next subset
$S_{i+1}$ in a subset chain~\eqref{chain}. The exceptional case occurs when
$1\not\in S_i$ and $S_{i+1}=\{1_G\}$, and, as we mentioned at the end of the
previous section, in this exceptional case we would typically use an
exhaustive search through $S_i$ to find the required `sifting element'. 

Our general approach assumes that we are able to test membership in each of the $S_i$ and to select a uniformly distributed 
random element from some subset `related to' $S_i$ in the chain~\eqref{chain}; see~Section~\ref{example}
for examples.

\begin{definition}\label{memtest}
A $4$-tuple $(G_0,H,K,\ismember)$ 
is said to satisfy the {\em membership test condition in $G$} 
if the following hold:
\begin{itemize}
\item[(a)] $G$ and $G_0$ are finite groups such that $G_0\leq G$;
\item[(b)] $H$ and $K$ are non-empty subsets of $G_0$ such that 
$H\supset K$.
\item[(c)] $\ismember$ is a 
one-sided Monte Carlo algorithm whose input is a pair $(y,e)$, where $y\in G$ and $e$ is a non-negative real number. It satisfies the following condition, either for all inputs $(y,0)$ (in which case it is a deterministic algorithm), or for all inputs $(y,e)$ with $0<e<1/2$. \em{The output is either {\sc True} or {\sc False}, and moreover, if $y\in K$ then the output is {\sc True}, and also
		$\prob{\mbox{output is {\sc True}}\ |\,y\in
H\setminus K}\leq e$.}
\end{itemize} 
\end{definition}

\noindent \textbf{Note:} For an enhanced version of an $\ismember$ test
giving back additional information for later use consult the examples
for $M_{11}$ and $\Ly$ in Section \ref{sect:appns}.

\smallskip
We show that if a $4$-tuple $(G_0,H,K,\ismember)$ satisfies the membership
test condition in a group 
$G$, then we can design an algorithm $\basicsift$ such that
$(G_0,H,K,\basicsift)$ satisfies the basic sift condition in $G$.
As mentioned above, we assume that we can select uniformly distributed random elements from some subset $L$ of $G$ `related to' the subset $H$. The most general conditions that the subset $L$ must satisfy are given in the following definition. 

\begin{definition}\label{def:hl}
Suppose that $G$ is a finite group and $H,\ K,\ L\subseteq G$.
We say that $(H,K,L)$ is a {\em sifting triple} if
\begin{equation}\label{eq:hl}
HL\subseteq H, \ \mbox{and, for all}\ h\in H,\
hL\cap K\neq\emptyset.
\end{equation}
\end{definition}

The reason why we introduce the subset $L$ in a sifting triple
is that it is rarely possible
to make random selections from arbitrary subsets of $G$, such as $H$,
but we can often make random selections from subgroups. Thus one
choice for $L$ is a subgroup satisfying (\ref{eq:hl}). Moreover we can
sometimes obtain a more efficient algorithm by restricting to a `nice
subset' $L$ of such a subgroup, provided that we can still make random
selections from $L$. Sometimes this is possible simply because $L$ is
small enough to hold in the memory. In that latter case we do not have to
perform a random search, but can use an exhaustive search. This is
analysed in Section~\ref{subsec:cosetreps}.

If $(H,K,L)$ is a sifting triple then the number 
$$
p(H,K,L)=\min_{h\in H}\frac{|hL\cap K|}{|L|}
$$
is called the {\em sifting parameter}. We note that the definition of a 
sifting triple implies that $p(H,K,L)>0$. The sifting parameter plays an
important r\^ole in estimating the complexity of Algorithm~\ref{sift}.

\subsection{A {\basicsift} algorithm using random search}\label{5.1}
\leavevmode

\addtocounter{algonr}{1}
\begin{algorithm}
\SetKw{Input}{Input:}
\SetKw{Output}{Output:}
\SetKw{Set}{set}
\SetKwFunction{IsMember}{\sc IsMember}
\SetKwFunction{Random}{\sc RandomElement}
\caption{A \basicsift\ algorithm using random search}
Algorithm~\ref{basicsift}: $\basicsiftrandom$\label{basicsift}
\tcc{See Theorem~\ref{prop:basic-mc} for notation}
\medskip
\Input{\rm $(x,\ve)$ where $x\in G$, and  $0< \ve<1/2$}\;
\Output{\rm $y$, where either $y=$ {\sc Fail}, or $y\in S$}\;
\Set $\displaystyle{e=\left\{\begin{array}{ll}
0&\mbox{if $\ismember$ is deterministic}\\
\ve p/(2(1-p))& \mbox{otherwise}
\end{array}\right.}$\;
\Set $\displaystyle{N=\left\{\begin{array}{ll}
\lceil\log(\ve)/\log (1-p)\,\rceil & \mbox{if $\ismember$ is deterministic}\\
\lceil\log(\ve/2)/\log (1-p)\,\rceil & \mbox{otherwise}\\
\end{array}\right.}$\;
\Set $n=0$\;
\Repeat{$n\geq N$}
	{
	\Set $y=\Random(L)$\;
	\If{\IsMember$(xy,e)$}{\Return $y$}
	\Set $n = n+1$}
\tcc{at this stage, none of the elements $y$ has been\\
returned during the for-loop}
\Return\fail
\end{algorithm}

\begin{theorem}\label{prop:basic-mc}
Suppose that $(G_0,H,K,\ismember)$ satisfies the membership test
condition in a group $G$ and that $L$ is a subgroup of 
$G_0$ such that $(H,K,L)$  is a sifting triple.
If $\randel(L)$ returns uniformly distributed, independent random
elements of $L$, and \basicsift\ is Algorithm~$\ref{basicsift}$, then
the $4$-tuple $(G_0,H,K,\basicsift)$
satisfies the basic sift condition in $G$. 
Moreover the cost of
executing $\basicsiftrandom(\cdot,\ve)$ is at most 
$$O\left(\log(\ve^{-1})\,
p^{-1}\left(\xi + \varrho+\nu(e)\right)\right), $$
where $p = p(H,K,L)$ and $\varrho,\ \xi$ and $\nu(e)$ are
upper bounds for the costs of a group operation in $G$, a random
selection from $L$, and one run of the procedure
$\mbox{\sc IsMember}(\cdot,e)$, respectively, where $e=0$ if $\ismember$
is deterministic, and $e=\ve p(H,K,L)/(2-2p(H,K,L))$ otherwise.
\end{theorem}
\begin{proof} 
If a group element $y$ is returned then, by (\ref{eq:hl}), $Hy\subseteq HL\subseteq H$. 

Let $E$ denote the event that ``the output of the procedure is either
\fail\ or an element $y$ with $xy\not\in K$''. We are required to show that
$\prob{E|x\in H}\leq\varepsilon$. 
Suppose that $x\in H$. For $i=0,\ldots,N-1$, let $E_i$ denote the event 
``the $(i+1)$-th execution of the procedure \randel\ occurs''; let $y_i$
denote the element $y$ returned by the $(i+1)$-th execution of \randel, and
let $z_i$ denote the result returned by the call to $\ismember(xy_i,e)$.
If $E_i$ does not occur for some $i$ then the values of $y_i$ and $z_i$ 
are undefined.
The event $E_i$ is the disjoint union of the following three events:
\begin{eqnarray*}
K_i&:& E_i\mbox{ and }xy_i\in K;\\
F_i&:& E_i\mbox{ and }xy_i\not\in K\mbox{ and }z_i=\false;\\
T_i&:& E_i\mbox{ and }xy_i\not\in K\mbox{ and }z_i=\true.
\end{eqnarray*}
Note that $E_i$ occurs if and only if, for each $j<i$, the event $E_j$ occurred
and $z_j=\false$, that is to say, $E_i=F_0\cap\cdots\cap F_{i-1}$. 
Similarly, given, $x\in H$, the event $E$ occurs if and only if either 
$F_0\cap F_1\cap\cdots\cap F_{N-1}$ occurs, or, for some $i$, 
each of $E_1,\ldots,E_i$ occurs, $xy_i\not\in K$ and $z_i=\true$.

Suppose now that $x\in H$, and let $y\in L$ such that $xy\not\in K$.
Then by (\ref{eq:hl}), $xy\in HL\subseteq H$, and hence $xy\in H\setminus K$.
By the definition of the membership test condition, the conditional
probability $e_0$ that the returned value of $\ismember(xy,e)$ is {\sc True}, given that $xy\in H\setminus K$, satisfies $0\leq e_0\leq e$.

Let $p$ denote the sifting parameter $p(H,K,L)$. 
Since we are making independent uniform random selections, 
we have, for each $i\leq N-1$, that the probability $\prob{K_i|E_i}$ is 
independent of $i$, 
and also that $$
\prob{K_i|E_i}\geq \frac{|xL\cap K|}{|L|}\geq p.
$$  
Set $p_0=\prob{K_i|E_i}$. Then, using the rule 
$\prob{A\cap B|C}=\prob{A|B\cap C}
\prob{B|C}$, 
$$
\prob{F_i|E_i}=\prob{xy_i\not\in K|E_i}\cdot\prob{z_i= \mbox{\sc False}\,|\,
E_i\mbox{ and }xy_i\not\in K}=(1-p_0)(1-e_0)
$$
with $e_0$ as defined above, and similarly $\prob{T_i|E_i}= (1-p_0)e_0$.

The procedure finishes when processing the $i$-th random element $y_i$ if it has not finished
while processing  $y_j$ for any $j<i$, and either $K_i$ or $T_i$ occurs.
In this situation, if $K_i$ occurs, then by the requirements of the
membership test condition, the procedure will return $y_i$ with $xy_i\in K$; similarly,
if $T_i$ occurs, then again the procedure will return $y_i$, but this time with $xy_i\not\in K$.
Thus the procedure returns the element $y_i$ with $xy_i\not\in K$ (for a particular value of $i$)
if and only if $F_0\cap\cdots\cap F_{i-1}\cap T_i=T_i$ occurs, and
$$
\prob{T_i}=e_0(1-p_0)((1-p_0)(1-e_0))^{i}. 
$$
It follows that the procedure returns an element $y\in L$ with $xy\not\in K$
if and only if $T_i$ 
occurs for some $i=0,\dots,N-1$, and the probability
of this is
$$
\sum_{i=0}^{N-1}e_0(1-p_0)^{i+1}(1-e_0)^{i}=e_0(1-p_0)\frac{1-(1-p_0)^N(1-e_0)^N}
{1-(1-p_0)(1-e_0)}
\leq \frac{(1-p_0)e_0}{p_0},
$$
since $1-(1-p_0)(1-e_0)=p_0+(1-p_0)e_0\geq p_0$.
Finally, the procedure returns {\sc Fail} if and only if the event $F_0\cap F_1\cap\dots
\cap F_{N-1}$ occurs and the probability of this is $(1-p_0)^N(1-e_0)^N$.

We derive the required estimates of these probabilities as follows.
Note that, since $p\leq p_0$ and $0\leq e_0\leq e$, we have
\[
\frac{(1-p_0)e_0}{p_0}=(p_0^{-1}-1)e_0\leq (p^{-1}-1)e
\]
and this is $0$ if $e=0$, and is $\ve/2$ otherwise.   Hence, the probability that the procedure returns an element $y\in L$, with $Hy\subseteq H$ and $xy\not\in K$, is $0$ if $\ismember$ is deterministic, and is at most $\ve/2$ otherwise. Similarly, the probability that the
procedure returns {\sc Fail} is 
\[
(1-p_0)^N(1-e_0)^N\leq (1-p)^N \leq \frac{\ve}{\delta},
\]
by the definition of $N$, where $\delta=1$ if $\ismember$ is deterministic, and $\delta=2$ otherwise. Thus $(G_0,H,K,\basicsift)$ satisfies the basic sift 
condition in~$G$.

Finally we estimate the cost.
For each run of the {\bf repeat} loop, 
first we select a random element of $L$ at a cost of at most 
$\xi$. Then 
we perform a group operation to compute
$xy$ and we run {\sc IsMember}$(xy,e)$ at a cost of at most $\varrho+\nu(e)$, where $e=0$ if $\ismember$ is deterministic, and  $e=\ve p/(2(1-p))$ otherwise.
The number of runs of the loop is at most $N$ and, by Lemma~\ref{rem:mc},
$N$ is $O(\log(\ve^{-1}) p^{-1})$.  Thus the upper bound for the cost is proved. (Note that, for $\ve<1/2$ we have that $\ve p/(2(1-p))<1/2$ also.) \qed
\end{proof}

As already explained before Theorem~\ref{prop:basic-mc}, we often work with
sifting triples $(H,K,L)$ in which $L$ is a subgroup of $G_0$. 
Usually, there will be another subgroup $L' < L$, which is
used to define $K$ and we have $KL' \subseteq K$. In this situation
the following concept applies.

\begin{definition}
Suppose that $G$ is a finite group and that $L,\ L'$ are subgroups of $G$. 
A non-empty 
subset $S$ of $L$ is said to be {\em left $L'$-uniform} if
$S$ has the same number of elements in each of the left $L'$-cosets in $L$.
In other words, $|S\cap \ell L'|$ is constant for
all $\ell \in L$. 
\end{definition}

A left $L'$-uniform subset in $L$ must contain a left
transversal for $L'$ in $L$.
Notice that $L$ is left $\{1_G\}$-uniform, and more generally, if 
$L'$ is a subgroup, then any left transversal for $L'$ in $L$ is
left $L'$-uniform. As will become clear in the next lemma, $L'$-uniform
sets $S$ have `nice' properties with respect to the calculation
of probabilities. 
In certain cases we need to consider 
sifting triples $(H,K,S)$ in which $S$ is a 
left $L'$-uniform subset in some subgroup $L$ for which $(H,K,L)$
 is also a sifting triple. 
We
show that in such cases
the sifting parameter 
$p(H,K,S)$ is independent of the subgroup $L'$ and the left
$L'$-uniform subset $S$, and depends only on the subgroup $L$.

\begin{lemma}\label{lem:p}
Let $(H, K, L)$ be a sifting triple 
in which $L$ is a subgroup,
let $L'$ be a subgroup of $L$ with $KL' \subseteq K$, and let $S$ be a  left $L'$-uniform subset of $L$. Then $(H,K,S)$ is also a sifting triple and 
$p(H,K,S)=p(H,K,L)$. 
\end{lemma}

\begin{proof} 
Since $HL\subseteq H$ and $S\subseteq L$, 
it follows that $HS\subseteq H$. 
Let $h\in H$. We shall show that $|hS\cap K|/|S|= |hL\cap
K|/|L|$. The result will then follow. By (\ref{eq:hl}), $hL\cap
K\ne\emptyset$. Note that, since $L'$ is a subgroup of $L$, and since $S$
is left $L'$-uniform, it follows that $L=SL'$, and $LL'=L$. In addition, we
have $KL'=K$. Thus Lemma~\ref{dede} implies that 
$(hL\cap K)L'= hL\cap K$, and in particular, $hL\cap K$ is a union of $r$ left $L'$-cosets, for some $r>0$. Each of these cosets is contained in $hL=hSL'$ and hence is of the form $hsL'$ for some $s\in S$. Thus $hL\cap K= \bigcup_{i=1}^r hs_iL'$ for some $s_1,\dots,s_r\in S$.

Further, since $S$ is left $L'$-uniform, the size $q=|s_iL'\cap S|$ is independent of $i$. Moreover, for each $i\leq r$, $hs_iL'\cap hS=h(s_iL'\cap S)$, and since $hS\subseteq hL$ it follows that
\[
hS\cap K = (hL\cap K)\cap hS = \bigcup_{i=1}^r (hs_iL'\cap hS)
=  \bigcup_{i=1}^r h(s_iL'\cap S),
\]
and therefore $|hS\cap K|=rq$. On the other hand, $hL\cap K = 
\bigcup_{i=1}^r hs_iL'$ has size $r |L'|$. Since $S$ has exactly $q$
elements in each of the left $L'$-cosets in $L$, we have $|S|=q |L:L'|$, and hence
\[
\frac{|hL\cap K|}{|L|}= \frac{r|L'|}{|L|} = \frac{rq}{|S|}
=\frac{|hS\cap K|}{|S|}
\]
proving the claim.
\qed
\end{proof}

\subsection{A {\basicsift} algorithm using a stored transversal}
\label{subsec:cosetreps}
\leavevmode

We now turn to a second general approach to designing a $4$-tuple that
satisfies the basic sift condition. 
This algorithm is defined for the case when we have a sifting triple
$(H,K,L)$ and a subgroup $L'\leq L$ as in Lemma~\ref{lem:p}. 
Unlike Algorithm~\ref{basicsift}, where we
choose elements of $L$ at random, Algorithm~\ref{basicsiftcr} 
deterministically tests every
element of a complete set $S$ of left coset representatives calculated
beforehand.
Thereby we turn the random search above into a deterministic exhaustive
search. As will be explained below, this can reduce the expected value
of the runtime significantly.

We use Algorithm~\ref{basicsiftcr} 
when the index of $L'$ in $L$, and thus the size of
$S$, is small enough to allow $S$ to be stored completely.
We still allow the use of randomised or deterministic
{\ismember} methods. In the latter case, the whole basic sift procedure
is deterministic.

We would like to draw attention to 
a little trick we use to simplify the analysis
of the error probability of Algorithm~\ref{basicsiftcr}.  We artificially 
introduce a randomly chosen order in which the coset representatives are
tried. This makes the analysis less dependent on the input group element.

\addtocounter{algonr}{1}
\begin{algorithm}
\SetKw{Input}{Input:}
\SetKw{Output}{Output:}
\SetKw{Set}{set}
\SetKwFunction{Random}{\sc RandomElement}
\SetKwFunction{IsMember}{\sc IsMember}
\caption{A \basicsift\ algorithm using a left transversal of $L'$ in $L$}
Algorithm~\ref{basicsiftcr}: {\basicsiftcosetreps}
\label{basicsiftcr}
\tcc{See Theorem~\ref{prop:basic-cr} for notation}
\medskip
\Input{\rm $(g,\ve)$ where $g\in G$, and $0 \leq \ve<1/2$}\;
\tcc{ $\ve=0$ if and only if {\ismember} is deterministic} 
\Output{\rm $y$, where either $y=$ {\sc Fail}, or $y\in S$}\;
\Set $\displaystyle{e=\left\{\begin{array}{ll}
0&\mbox{if $\ismember$ is deterministic}\\
\min\left\{\displaystyle{\ve\cdot\frac{n+1}{k-n}},\frac{1}{3}\right\} & 
\mbox{otherwise, where } k = |S|,\ 
n = \min_{h \in H} |hS\cap K|
\end{array}\right.}$\;
\Set $T = S$\;
\For {$i=1,2,\ldots,k$}
        {
        \Set $y=\Random(T)$\;
	\If{\IsMember$(gy,e)$}{\Return $y$}
        \Set $T = T \setminus \{y\}$\;
        }
\tcc{ we only reach this stage if $g \notin H$, because otherwise one of
the {\ismember} tests must have returned {\sc True}}
\Return\fail
\end{algorithm}

\begin{theorem}\label{prop:basic-cr}
Suppose that $(G_0,H,K,\ismember)$ satisfies the membership test condition 
in a group $G$.
Assume further that $L$ is a subgroup of $G_0$, 
such that $(H,K,L)$ is a sifting 
triple, that $L' < L$ with $KL' = K$, and that
$S = \{s_1, \ldots, s_k\}$ is a left 
transversal of $L'$ in~$L$. If,  for any $T \subseteq S$,  $\randel(T)$ returns 
uniformly distributed, independent random elements of $T$,
and {\basicsift} is Algorithm~$\ref{basicsiftcr}$, then the $4$-tuple
$(G_0,H,K,\basicsift)$ satisfies the basic sift condition in $G$.

The cost of executing $\basicsiftcosetreps(\cdot,\ve)$ is less than
$k \cdot (\xi_S+\varrho + \nu(e))$ where $\xi_S$ is an upper bound for
the cost of selecting a random element from 
a subset of $S$, $\varrho$ and $\nu(e)$ are upper bounds
for the costs of a group operation in $G$, and one run of the procedure
$\ismember(\cdot,e)$, respectively. Here $e=0$ if $\ismember$ is
deterministic, and $e = \min\left\{\ve(n+1)/(k-n),1/3\right\}$ 
otherwise, where $n = \min_{h \in H} |hS \cap K|$.
\end{theorem}
\begin{proof}
We remark first, that for every $g \in H$ there is an element $l \in L$
such that $gl \in K$ by hypothesis~(\ref{eq:hl}). As $S$ is a left
transversal for $L'$ in $L$, there are 
$s \in S$ and $l' \in L'$ such that $l=sl'$. Now $gsl'\in K$, and so $gs\in
KL'=K$. 
Therefore, if $g \in H$, then Algorithm~\ref{basicsiftcr} cannot return {\sc Fail}, as the {\ismember} test is
one-sided Monte Carlo. Also, this argument proves all statements
in the theorem in the case where {\ismember} is deterministic. 

Thus from now on we will assume that {\ismember} is not deterministic, and
therefore that $0 < \varepsilon < 1/2$, and hence $e$ is non-zero.

As $HL = H$, the set $H$ is a union of left $L$-cosets,
and, a fortiori, also a union of left $L'$-cosets. Analogously, 
$KL' = K$ means that $K$ is a union of left $L'$-cosets, and, of course, so is $gL\cap K$.
 For any
given $g$, the algorithm looks for a random 
element $y$ in $S \subset L$ such that $gy\in K$; in other words, 
it searches the coset $gL$ for elements of $K$. Thus, the number of
elements $s \in S$ with $gs \in K$ is equal to the number of
left $L'$-cosets contained in $gL \cap K$.
Let $g\in H$.
As, by Lemma~\ref{lem:p}, $p(H,K,S)=p(H,K,L)$, 
and $|gL\cap K|/|L| = |gL\cap K|/(k|L'|)$
we obtain that  
$$
|gS \cap K|=
|S|\frac{|gS\cap K|}{|S|}\geq |S|\min_{h\in H}
\frac{|hS\cap K|}{|S|}=kp(H,K,S)=kp(H,K,L).
$$

Let $E$ denote the event ``the procedure returns $y\in S$ with $gy\not\in K$''.
To check the basic sift condition for Algorithm~\ref{basicsiftcr}
in the case of a randomised {\ismember} test, we have to show that $\Prob(E
\mid  g\in H) \leq \varepsilon$.

Suppose now that  $g\in H$.
For $i=1,\ldots,k$, let $E_i$ denote the event: ``the $i$-th  execution of the
procedure {\sc RandomElement} occurs''; 
let $y_i$ denote
the element $y$ returned by the $i$-th execution of the procedure
$\randel$, and let $z_i$  denote the result returned
by the call to $\ismember(gy_i,e)$ (for steps $i$ that do not happen, $y_i$
and $z_i$ are undefined). 

Then $E_i$ is the disjoint union of the following three events:
\begin{eqnarray*}
  K_i & : & E_i \mbox{ and } gy_i\in K; \\
  F_i & : & E_i \mbox{ and } gy_i\not\in K \mbox{ and } z_i = \false; \\
  T_i & : & E_i \mbox{ and } gy_i\not\in K \mbox{ and } z_i = \true.
\end{eqnarray*}

Note that $E_i$ occurs if and only if, for each $j<i$, the event $E_j$
occurred and $z_j=\false$, that is to say, $E_i = F_1\cap \cdots \cap
F_{i-1}$. Similarly, given $g\in H$, the event $E$ occurs if and only if,
for some $i$, each of $E_1,\ldots, E_i$ occurs, $gy_i\not\in K$, and 
$z_i=\true$. Thus,
given $g\in H$, the event $E$ occurs if and only if, $F_1\cap\cdots\cap
F_{i-1}\cap T_i=T_i$ occurs for some $i$ with $1 \leq i \leq k$.

Since in step $i$ we choose $y_i$ only among
those coset representatives that have not been tried before and
we only reach step $i$ if $gy_j \notin K$ for $1 \leq j < i$, 
the probability $\prob{gy_i \notin K \mid  E_i}$ is not
independent of $i$.
Namely, \[ \prob{gy_i \notin K \mid  E_i } = (k+1-i-n_g)/(k+1-i) \] where 
$n_g = |gS \cap K|$, as in step $i$ there are $k+1-i$ coset
representatives in the set $T$ of which $k+1-i-n_g$ do not multiply
$g$ into $K$. 

It is easy to see that
\[
  \prob{F_i\mid E_i} = \prob{gy_i \notin K\mid E_i} \cdot 
  \prob{z_i = \false \mid gy_i \notin K \mbox{ and } E_i }
  \]
and 
so
\[  
  \prob{F_i \mid E_i} \leq \frac{k+1-i-n_g}{k+1-i}. \]
Similarly we have
\[
\prob{T_i \mid E_i} \leq \frac{k+1-i-n_g}{k+1-i} \cdot e.
\]

As in the proof of Theorem~\ref{prop:basic-mc}, Algorithm~\ref{basicsiftcr}
finishes in step $i$, if it has not finished in an earlier step, and
 $K_i$ or $T_i$ occurs. In this situation, if $K_i$
occurs, then the procedure will return $y_i$ with $gy_i \in K$, which is
a correct result. Therefore, an error produced by step $i$ occurs exactly
in the event $T_i$, and
\[
  \prob{T_i}
  \leq \left( \prod_{j=1}^i \frac{k+1-j-n_g}{k+1-j} \right) \cdot e .
\]
Moreover, 
no error can possibly occur in step $i$ for $i > k-n_g$.

Therefore, for an input $(g,\ve)$ with $g\in H$, the total probability
that Algorithm~\ref{basicsiftcr} returns an element $y \in S$ with $gy
\notin K$ is
\[
  \sum_{i=1}^{k-n_g} \left( \prod_{j=1}^i \frac{k+1-j-n_g}{k+1-j} \right)
  \cdot e .
\]
Note 
that, for $i=1,\ldots,k-n_g$,
$$
\prod_{j=1}^i \frac{k+1-j-n_g}{k+1-j}=\frac{(k-n_g)(k-n_g-1)\cdots(k-n_g-i+1)}
{k(k-1)\cdots(k-i+1)}=\frac{(k-i)!}{k!}\cdot\frac{(k-n_g)!}{(k-i-n_g)!}.
$$
Hence
$$
  \sum_{i=1}^{k-n_g} \left( \prod_{j=1}^i \frac{k+1-j-n_g}{k+1-j} \right)
 =\sum_{i=1}^{k-n_g} \frac{(k-i)!}{k!} \cdot \frac{(k-n_g)!}{(k-i-n_g)!}
                     \cdot \frac{n_g!}{n_g!}
 = \binom{k}{n_g}^{-1} \cdot \sum_{i=1}^{k-n_g} \binom{k-i}{n_g}.
$$
We can
simplify the sum further by repeated use of the well known summation
formula for binomial coefficients:
$$
  \binom{a}{b} + \binom{a}{b-1} = \binom{a+1}{b}.
$$
The last summand (with $i=k-n_g$)  is equal to $\binom{n_g}{n_g} = 1=
\binom{n_g+1}{n_g+1}$. In the latter form it can be added to
the second last summand resulting in $\binom{n_g+2}{n_g+1}$. This
can be repeated until the first summand, thereby proving that 
$$
\sum_{i=1}^{k-n_g} \binom{k-i}{n_g}=\binom{k}{n_g+1}.
$$
This, however, implies that the total
probability of an error is
\[
  \binom{k}{n_g}^{-1} \cdot \binom{k}{n_g+1} \cdot e
  = \frac{n_g! \cdot (k-n_g)!}{k!} \cdot \frac{k!}{(n_g+1)! \cdot (k-n_g-1)!}
    \cdot e
  = \frac{k-n_g}{n_g+1} \cdot e.
\]
Thus, as $n\leq n_g$, 
for an arbitrary element $g \in H$, 
the error probability is bounded by
$$
  \frac{k-n}{n+1} \cdot e \leq \ve.
$$

As for the cost, the loop terminates at the latest after $k$ steps, each of 
which
has a random element selection from $T$, 
one group multiplication for computing $gy_i$, and one call to $\ismember$.
\qed
\end{proof}

Our hypotheses in Theorem~\ref{prop:basic-cr}
imply that $S$ is $L'$-uniform. However, since
we want to store $S$ completely, there is no point in choosing left
$L'$-uniform sets with two or more elements in each left $L'$-coset
of $L$.

\subsection{Comments on and comparison of Algorithms~\ref{basicsift} 
and \ref{basicsiftcr}} \leavevmode

To compare Algorithms~\ref{basicsift} and~\ref{basicsiftcr}, assume 
that $L$ is a subgroup and 
we want to sift from a set $H$ with $HL=H$ down to a set $K$ with
$KL'=K$, and that $L' < L$ with $[L:L'] = k$. Then we can either  use 
Algorithm~\ref{basicsift} 
or use Algorithm~\ref{basicsiftcr} with $S$
being a left transversal of $L'$ in $L$. Recall 
that $p(H,K,L)=p(H,K,S)=p$, say (see Lemma~\ref{lem:p}). 
Let $k$ denote the index $|L:L'|$, and
let $n$ denote $\min_{h\in H}|hS\cap K|=pk$.
In the second 
case we have to calculate and store $S$ beforehand.
In Algorithm~\ref{basicsiftcr}, once we compute that 
a random element $y$
does not 
multiply $g$ into $K$, $y$ cannot be selected again by a subsequent
call of $\randel$.
Therefore we expect that Algorithm~\ref{basicsiftcr} performs
better than Algorithm~\ref{basicsift} in this situation.

In Algorithm~\ref{basicsift} the bound for the error probability
in all calls of the {\ismember} test is
$ e_1 = \ve p/(2-2p) = \ve n/(2(k-n))$
(recall that $p = n/k$), 
whereas in Algorithm~\ref{basicsiftcr}
the bound for the error probability for the $\ismember$ calls is
$e_2 = \ve(n+1)/(k-n)$ (at least when $\ve$ is not too big so that
$e_2$ is not defined to be $1/3$), which is a little bit more than 
$2e_1$.
Thus, due to the deterministic nature of the choice of $y$ in
Algorithm~\ref{basicsiftcr}, we can afford bigger error bounds for the
$\ismember$ tests.
Further, the expected number of steps in Algorithm~\ref{basicsift}
is $1/p$ (geometric distribution), which is $k/n$ as $p=n/k$. The
expected number of steps in Algorithm~\ref{basicsiftcr}
is $(k+1)/(n+1)$. 

These calculations suggest that, whenever
it is possible to store all elements of $S$,
Algorithm~\ref{basicsiftcr} should be preferred over
Algorithm~\ref{basicsift}.

If the $\ismember$ test is deterministic and happens to work not only for
elements of $H$, but also for arbitrary elements of $HS$, then one can
dispense with the hypothesis $HS \subseteq H$ altogether and apply
Algorithm~\ref{basicsiftcr} verbatim for
any set $S \subseteq G$ satisfying $hS \cap K \not= \emptyset$ for
all $h \in H$. In this case Algorithm~\ref{basicsiftcr} will be a fully
deterministic algorithm with guaranteed finite runtime of at most $|S|$
steps.

\section{{\sc BasicSift}: with special subsets $H$ and $K$}\label{sect:CL}

In this section we describe a rather general situation where 
the conditions in 
(\ref{eq:hl}) are guaranteed to hold. The conditions on the subsets $H$, $K$ of 
the finite group $G$ are as follows:
\begin{equation}\label{localhk}
H=CL=C'L,\ K=C'L',\ \mbox{where}\ L'<L\leq G,\   \mbox{and}\ C,\ C'\subseteq G, \ \mbox{with}\ C,\ C'\ne\emptyset.
\end{equation}

Under these conditions we derive also a new expression for the sifting
parameter $p(H,K,L)$ required for Algorithm~\ref{basicsift} and Theorem~\ref{prop:basic-mc}.

\begin{proposition}\label{prop:localhk}
Let $G$, $L$, $L'$, $C$, $C'$, $H$, and $K$ be as above so that~{\rm(\ref{localhk})} holds. Then $K\subseteq H$, and if $H\neq K$ then $(H,K,L)$ is a sifting
triple. Further, 
$$
p(H,K,L)=\min_{y\in\mathcal Y}\frac{|yL\cap K|}{|L|}=\min_{y\in\mathcal Y}\frac{|(yL\cap C')L'|}{|L|},
$$
where $\mathcal Y$ is a set of representatives in $C'$ for the left
$L$-cosets contained in $H$.
\end{proposition} 

\begin{proof}
Since $L'\subseteq L$ we have $C'L'\subseteq C'L=CL$, that is, $K\subseteq H$.
Note that, since $1\in L'$, we have $C'\subseteq K$ and $C\subseteq H$.

Suppose now that $K\neq H$. Since $H=CL$ and $L$ is a subgroup, it follows that $HL\subseteq H$. Let $y\in H$. To complete the proof of~\eqref{eq:hl}, we need to show that $yL\cap K$ is non-empty.
Since $y\in H$ and $H=CL=C'L$ we have $y=ck$ where $c\in C'$, $k\in L$, and hence $c=yk^{-1}$ and $c\in yL\cap C'$. As $C'\subseteq K$, we obtain
$c\in yL\cap K$. Thus $yL\cap K\ne\emptyset$.

Now it only remains to show that the assertion in the displayed line of
the proposition is valid. It follows from~\eqref{localhk} that, for $y\in H$, 
$yL\cap K=yL\cap (C'L')=(yL\cap C')L'$,
by Dedekind's modular law (Lemma~\ref{dede}). Hence, for all $y\in H$, we have
$$
\frac{|yL\cap K|}{|L|}=\frac{|(yL\cap C')L'|}{|L|}.
$$
Suppose that $y\in H$ and $y=ck$ where $c\in C'$ and $k\in L$. Then 
$yL\cap K=cL\cap K$ and so the minimum value of
$|yL\cap K|/|L|$ over all $y\in H$ is equal to the minimum value of $|cL\cap K|/|L|$ over all $c\in \mathcal{Y}$.
The displayed assertion follows.
\qed
\end{proof}

We will apply Algorithm~\ref{basicsift} with $H$, $K$ 
as in (\ref{localhk}) in the following context: $G_0$ is a subgroup of a 
finite group $G$, 
the group $G_0$ has a 
descending subgroup chain
\begin{equation}\label{Lchain}
G_0=L_0 > L_1 > \dots > L_{k}=\{1\},
\end{equation}
and also has a sequence of non-empty subsets 
\begin{equation}\label{Cchain}
C_0=\{1\},C_1,\ldots,C_k\quad\mbox{such that}\quad C_{i+1}L_i=C_iL_i
\quad\mbox{for all}\quad i<k. 
\end{equation}
Thus 
(\ref{localhk}) holds for $(H,K)=(C_iL_i,C_{i+1}L_{i+1})$ for each $i<k$. By 
Proposition~\ref{prop:localhk}, we have a descending chain
\begin{equation}\label{CLchain}
G_0=C_0L_0\supseteq C_1L_1\supseteq\dots\supseteq C_{k}L_{k}=C_{k}
\end{equation}
and by Proposition~\ref{prop:localhk},  Algorithm~\ref{basicsift} applies 
to each of the pairs $(C_iL_i,C_{i+1}L_{i+1})$  such 
that $C_iL_i\ne C_{i+1}L_{i+1}$ ($0\leq i<k$). Thus if, for
$i=0,\ldots,k-1$, the $4$-tuple
$(G_0,C_iL_i,C_{i+1}L_{i+1},\ismember_i)$ satisfies the membership test 
condition in $G$
for some algorithm $\ismember_i$, and if we have an algorithm
$\basicsift_{k}$ such that the $4$-tuple $(G_0,C_{k},\{1\},\basicsift_{k})$ satisfies the
basic sift condition in $G$, then we may use the procedures $\basicsift_i$ 
in Algorithm~\ref{sift}.
If $|C_k|$ is small, $\basicsift_k$ may simply test each member of
$C_k$ for equality with the input element (if $1\not\in C_k$), or its
inverse (if $1\in C_k$).
 The next two sections offer some possibilities for these
procedures that have been effective for computing with some of the sporadic
simple groups.

\section{{\sc IsMember} using conjugates}\label{sect:conj}

In this section we apply the theory developed in Sections~\ref{sect:mc} and~\ref{sect:CL}, especially in Section~\ref{sect:CL}, to sift an element down a subgroup chain such as (\ref{Lchain})
making use of an auxiliary subset sequence. 
This application uses conjugates of
an element $a$ with the following property:
\begin{equation}\label{a}
\begin{array}{l}
a\in L_{k-1}\setminus\{1\}\ \mbox{such that, for each $i=0,\dots,k-2$,}\\ 
\mbox{each $L_i$-conjugacy class in $a^{G_0}\cap L_i$ intersects $L_{i+1}$ 
non-trivially.}
\end{array}
\end{equation}
We construct an associated subset sequence (\ref{Cchain}) recursively as follows.
The first subset is $C_0=C_{G_0}(a)\mathcal T_0$ where $\mathcal{T}_0=\{1
\}$. Consider a typical link in the chain (\ref{Lchain}), say $L_i>L_{i+1}$ for
$i\leq k-2$, and 
suppose that 
we have already constructed  the subset $C_i$ corresponding to $L_i$, and $C_i$ is 
of the form $C_i=C_{G_0}(a)\mathcal{T}_i$, where $\{a^y\,|\,y\in\mathcal{T}_i\}$ is a 
set of $L_i$-conjugacy class representatives in $a^{G_0}\cap L_i$. 
Then 
$a^{G_0}\cap L_{i+1}=\bigcup_{y\in\mathcal{T}_i}(a^{yL_i}\cap L_{i+1})$, and by condition (\ref{a}),
each $a^{yL_i}\cap L_{i+1}$ is non-empty. For each $y\in\mathcal{T}_i$, choose 
$\mathcal{U}(y)\subset L_i$ such that $\{a^{yu}\,|\,u
\in\mathcal{U}(y)\}$ is a set of representatives for the $L_{i+1}$-conjugacy 
classes in $a^{yL_i}\cap L_{i+1}$.
Define $\mathcal{T}_{i+1}=\bigcup_{y\in\mathcal{T}_i}\, y\mathcal{U}(y)$, and define
the subset $C_{i+1}$ corresponding to $L_{i+1}$ by $C_{i+1}=C_{G_0}(a)\mathcal{T}_{i+1}$. In addition set $C_k=\{1\}$. 

We prove that (\ref{localhk}) holds, and we also derive two expressions for the sifting parameter $p(H,K,L)$ required for Algorithm~\ref{basicsift} and Theorem~\ref{prop:basic-mc}. The first expression shows that $p(H,K,L)$ is a ratio of the sizes of two
special subsets of conjugates of the element $a$, while the second expression provides a means of computing $p(H,K,L)$ from the orders of various centraliser subgroups.

\begin{proposition}\label{prop:conj} 
Suppose that $G$, $G_0$, $a$, 
$L_i$, $L_{i+1}$, $C_i$, $C_{i+1}$, $\mathcal T_i$, 
$\mathcal T_{i+1}$, and the $\mathcal{U}(x)$, for $x\in\mathcal{T}_i$, are as at the beginning of this section, and set $H=C_iL_i$, $K=C_{i+1}L_{i+1}$, $L=L_i, L'=L_{i+1}$, $C=C_i$ and $C'=C_{i+1}$. Then $\mathcal T_{i+1}L_i=\mathcal T_iL_i$ and~\eqref{localhk} holds, and if also $H\ne K$, then $(H,K,L)$ is a sifting
triple. 
Further, 
$$
p(H,K,L)=\min_{x\in\mathcal T_i}\frac{|a^{xL_i}\cap L_{i+1}|}{|a^{xL_i}|}=\frac{1}{|L:L'|}\min_{x\in\mathcal T_i}
\Big\{ |C_L(a^x)| \sum_{u\in\mathcal{U}(x)}\frac{1}{|C_{L'}(a^{xu})|}\Big\}.
$$
\end{proposition}

\begin{proof}  
By the definition of $\mathcal{T}_{i+1}$, we have that 
$\mathcal{T}_{i+1}\subseteq \mathcal{T}_iL_i$. 
Also, since~\eqref{a} holds, for each $x\in\mathcal T_i$
there exists $k\in L_i$ such that $xk\in\mathcal T_{i+1}$. 
Thus $\mathcal
T_i\subseteq\mathcal T_{i+1}L_i$, and so, since $L_i$ is a subgroup, we have 
$$
\mathcal T_iL_i\subseteq (\mathcal T_{i+1}L_i)L_i=\mathcal T_{i+1}L_i\subseteq (\mathcal
T_iL_i)L_i=\mathcal T_iL_i.
$$ 
Hence $\mathcal T_{i+1}L_i=\mathcal T_iL_i$. 
To prove~\eqref{localhk} it is sufficient to prove that $H=C'L=C_{i+1}L_i$. From the definition of $H$ 
we have 
\[
H=C_iL_i=C_{G_0}(a)\mathcal{T}_iL_i = 
C_{G_0}(a)\mathcal{T}_{i+1}L_i = C_{i+1}L_{i}=C'L.
\]
Thus~\eqref{localhk} holds. Moreover, if $H\neq K$, then,
by Proposition~\ref{prop:localhk}, then $(H,K,L)$ is a sifting triple.

It remains to show that the value of the sifting parameter $p(H,K,L)$ is as claimed.
Suppose that $h\in H$, and that $h=cxk$ with $c\in C_{G_0}(a)$, $x\in\mathcal 
T_i$, and $k\in L_i$. We claim that $|hL_i\cap K|=|(xL_i\cap C_{i+1})L_{i+1}|$. As $k\in L_i$, 
we certainly have $hL_i\cap K=cxL_i\cap K$. An easy calculation shows
that $cxL_i\cap C_{G_0}(a)\mathcal T_{i+1}L_{i+1}=c(xL_i\cap C_{G_0}(a)\mathcal T_{i+1}L_{i+1})$, and so $|cxL_i\cap C_{i+1}L_{i+1}|=|xL_i\cap C_{i+1}L_{i+1}|$. Therefore $|hL_i\cap K|=|xL_i\cap K|$. Finally, by Dedekind's modular law (Lemma~\ref{dede}, which applies since $(xL_i)L_{i+1}\subseteq xL_i$), we obtain
\[
xL_i\cap K=xL_i\cap C_{i+1}L_{i+1}=(xL_i\cap C_{i+1})L_{i+1}
\]
proving our claim.

Next we show that 
$xL_i\cap C_{i+1}=xC_{L_i}(a^x)\mathcal U(x)$, with $\mathcal{U}(x)$ as defined before Proposition~\ref{prop:conj} (recall that $x\in\mathcal{T}_i$).
Let $y\in xL_i\cap C_{i+1}$, so that $y=xk$ for some $k\in L_i$ and  $xk\in C_{i+1}$. 
Since $C_{i+1}= C_{G_0}(a)\mathcal T_{i+1}$, it follows that 
$a^{xk}\in a^{\mathcal T_{i+1}}\cap a^{xL_i}$. By the definition of 
$\mathcal T_{i+1}$, there is some $u\in \mathcal U(x)$ such that
$a^{xk}= a^{xu}$, and so $k\in C_{L_i}(a^x)u$. Therefore $xk\in xC_{L_i}(a^x)u$, and we obtain that $y=xk\in xC_{L_i}(a^x)\mathcal U(x)$. 
Conversely consider $y=xcu$, where 
$c\in C_{L_i}(a^x)$ and $u\in\mathcal U(x)$. 
As $\mathcal U(x)\subseteq L_i$,  
we have $y=xcu\in xL_i$. 
Further, $a^{xcu}=a^{xu}\in a^{x\mathcal U(x)}\subseteq a^{\mathcal T_{i+1}}$. 
Thus $y=xcu\in C_{G_0}(a)\mathcal T_{i+1}=C_{i+1}$. Therefore our claim is proved.

Putting the calculations in the last two paragraphs together, we have shown, for $h=cxk$ with $c\in C_{G_0}(a)$, $x\in\mathcal T_{i}$, and $k\in L_i$, that $|hL_i\cap K|= |xC_{L_i}(a^x)\mathcal U(x)L_{i+1}|.$
Now we calculate the size of $xC_{L_i}(a^x)\mathcal U(x)L_{i+1}$. We first observe that $xC_{L_i}(a^x)\mathcal U(x)L_{i+1}$
is a union of left $L_{i+1}$-cosets, and hence, it suffices 
to compute the number of such cosets contained in $xC_{L_i}(a^x)\mathcal U(x)L_{i+1}$. If $u_1$ and $u_2$ are distinct elements of $\mathcal U(x)$, then
$a^{xC_{L_i}(a^x)u_1L_{i+1}}=a^{xu_1L_{i+1}}$ and $a^{xC_{L_i}(a^x)u_2L_{i+1}}=a^{xu_2L_{i+1}}$, and so it follows from the definition of $\mathcal U(x)$ that
$a^{xC_{L_i}(a^x)u_1L_{i+1}}$ and $a^{xC_{L_i}(a^x)u_2L_{i+1}}$ are distinct conjugacy
classes in $L_{i+1}$. Thus $xC_{L_i}(a^x)u_1L_{i+1}$ and $xC_{L_i}(a^x)u_2L_{i+1}$ are
disjoint. 
Therefore $xC_{L_i}(a^x)\mathcal U(x)L_{i+1}$ is the disjoint 
union, over all $u\in\mathcal{U}(x)$, of $xC_{L_i}(a^x)uL_{i+1}$. 
Let $c_1, c_2\in C_{L_i}(a^x)$. Then 
$xc_1uL_{i+1}=xc_2uL_{i+1}$ if and only if $c_2^{-1}c_1\in uL_{i+1}u^{-1}$. Thus the number 
of left $L_{i+1}$-cosets in   $xC_{L_i}(a^x)uL_{i+1}$ is 
$|C_{L_i}(a^x)|/|C_{uL_{i+1}u^{-1}}(a^x)|
= |C_{L_i}(a^x)|/|C_{L_{i+1}}(a^{xu})|$. 
Hence, the definition of $\mathcal U(x)$ implies that
\begin{eqnarray*}
|hL_i\cap K|&=&|xC_{L_i}(a^x)\mathcal U(x)L_{i+1}|=\sum_{u\in \mathcal U(x)}|xC_{L_i}(a^x)uL_{i+1}|\\
&=&\sum_{u\in \mathcal U(x)}\frac{|C_{L_i}(a^x)|\cdot|L_{i+1}|}{|C_{L_{i+1}}(a^{xu})|}
=|C_{L_i}(a^x)|\sum_{u\in \mathcal U(x)}
\frac{|L_{i+1}|}{|C_{L_{i+1}}(a^{xu})|}\\
&=&|C_{L_i}(a^{x})|\sum_{u\in \mathcal U(x)}|(a^{xu})^{L_{i+1}}|
=|C_{L_i}(a^{x})|\cdot|a^{xL_i}\cap L_{i+1}|.
\end{eqnarray*}
Thus $$
\frac{|hL_i\cap K|}{|L_i|}=\frac{|C_{L_i}(a^{x})|\cdot|a^{xL_i}\cap L_{i+1}|}{|L_i|}=\frac{|a^{xL_i}\cap L_{i+1}|}{|a^{xL_i}|}
$$
and also
$$
\frac{|hL_i\cap K|}{|L_i|}=\frac{|C_{L_i}(a^{x})|}{|L_i:L_{i+1}|}
\sum_{u\in\mathcal{U}(x)}\frac{1}{|C_{L_{i+1}}(a^{xu})|}.
$$
Therefore we obtain that the displayed assertions for the sifting parameter 
also hold.
\qed
\end{proof}

The main benefit of working with conjugates is that, using the notation of Proposition~\ref{prop:conj},  
membership of $x$ 
in $H$ or $K$ is equivalent to membership of $a^x$ in $L_i$ 
or $L_{i+1}$, respectively; see Lemma~\ref{lem:conj1}.  
It is often easier to test whether a random conjugate of a known element
lies in a subgroup than to test membership of a random element in a subgroup or subset. This is true in particular if we have detailed 
information about subgroups of $L_i$ or $L_{i+1}$ generated by two $a$-conjugates.

\begin{algorithm}
\SetKw{Input}{Input:}
\SetKw{Output}{Output:}
\SetKwFunction{IsMember}{IsMember}
\caption{An \ismember\ algorithm for subsets}
\label{newismember}
Algorithm~\ref{newismember}: $\newismember$\\
\tcc{see Lemmas~\ref{lem:conj1} and~\ref{lem:conj2} for notation}
\medskip
\Input{\rm $(x,e)$ where $x\in G$, and $e=0$ if $\ismember$ is deterministic, and $0< e<1/2$ otherwise}\;
\Output{\rm \true\ or \false}\;
\Return \ismember$(a^x,e)$
\end{algorithm}

\begin{lemma}\label{lem:conj1}
Let $G$, $G_0$, $a$, $L_i$, $L_{i+1}$, $C_i$, $C_{i+1}$, $\mathcal T_i$, $\mathcal T_{i+1}$,  be as in  Proposition~{\rm\ref{prop:conj}}, set $H=C_iL_i$ and $K=C_{i+1}L_{i+1}$, and let $x\in G$. 
\begin{enumerate}
\item[(a)] 
The element $x\in H$ if and only if $a^x\in L_i$, and similarly, 
$x\in K$ if and only if 
$a^x\in L_{i+1}$.
\item[(b)] 
If $(G_0,L_i,L_{i+1},\ismember)$ satisfies the membership test condition in $G$, for
some 
algorithm $\ismember$, then 
so does $(G_0,H,K,\newismember)$ where the algorithm $\newismember$ is  given by \mbox{Algorithm}~$\ref{newismember}$.
\end{enumerate}
\end{lemma}

\begin{proof} It follows from the definition of $\mathcal{T}_i$ that $a^{
\mathcal{T}_iL_i}= a^{G_0}\cap L_i$. The first assertion in part (a) is then obvious, and
the second follows similarly. 

To prove part (b), recall the second assertion of part (a), namely that $x\in K$ if and only if
$a^x\in L_{i+1}$. If this condition holds then
the membership test condition (see Definition~\ref{memtest}) on $\ismember$ implies that $\ismember(a^x,e)=\true$ and hence we obtain $\newismember(x,e)=\true$.
Also, by part (a), $x\in H\setminus K$ if and only if
$a^x\in L_i\setminus L_{i+1}$. 
By the membership test condition on $\ismember$ we have 
$$
\prob{\mbox{output of $\ismember$ is }\true\,|\,x\in H\setminus K}\leq e
$$
and hence by the `definition' of $\newismember$ in Algorithm~\ref{newismember}, 
$$
\prob{\mbox{output of $\newismember$ is } \true\,|\,a^x\in L_i\setminus L_{i+1}}\leq e.
$$ 
Thus the membership test condition holds for $(G_0,H,K,\newismember)$ in $G$.
\qed \end{proof}

By Lemma~\ref{lem:conj1},  we can use $\newismember(a^{xy},e)$ to replace the algorithm $\ismember(xy,e)$ in the $\basicsift$ Algorithm~\ref{basicsift}. 
Some explicit instances of {\sc IsMember}
will be dis\-cussed in Section~\ref{sect:appns}. We discuss here one special 
case, namely where $L_{i+1}=\la a\ra$. Here it turns out that Lemma~\ref{lem:conj1} applies with $K=N_{G_0}(\la a\ra)$.
Before proving this assertion in Lemma~\ref{lem:conj2} below, we make a few
comments about the context in which it will arise.
(This context below occurs in several applications to sporadic simple groups.)

If condition (\ref{a}) holds for
a subgroup chain (\ref{Lchain}), then we construct, as at the beginning of 
this section, subsets $\mathcal{T}_i$ and $C_i=C_{G_0}(a)\mathcal{T}_i$, for each 
$i$, such that (\ref{Cchain}) and (\ref{CLchain}) both hold. Note that $a^{G_0}
\cap L_i=a^{\mathcal{T}_iL_i}$ and that  $\mathcal T_iL_i= \mathcal T_{i+1}L_i$ for each $i$; see  Proposition~\ref{prop:conj}.  Also $\la a\ra\leq L_{k-1}\leq L_i$, for all $i\leq k-1$. This means that $a^{G_0}\cap L_{i}$ contains $a$, and hence contains $a^{L_{i}}$. Thus ${\mathcal{T}}_{i}$ contains an element of $C_{G_0}(a)L_{i}$. In particular, if $L_{i}\leq C_{G_0}(a)$, then $\mathcal{T}_{i}$ contains an element of $C_{G_0}(a)$. (Note, however, that this element of $\mathcal{T}_i$ need not be equal to $1$.) 

It is tempting
to consider refining the chain (\ref{Lchain}) by inserting the subgroup
$\la a\ra$ to obtain a new chain with second last subgroup 
equal to $\la a\ra$. However condition (\ref{a}) may fail to hold for this new chain. 
For example if the original $L_{k-1}\cong \Z_2\times \Z_2$ then $a$ is an 
involution, and $|\mathcal{T}_{k-1}|=3$, but only one of the three
$L_{k-1}$-conjugacy classes in $a^{G_0}\cap L_{k-1}$ meets $\la a\ra$ non-trivially.
Nevertheless, the situation $L_{k-1}=\la a\ra$ arises often in applications, so
we end this section by extending the framework to include this case.

\begin{lemma}\label{lem:conj2}
Suppose that $G$, $G_0$, $a$, $L_i$, $L_{i+1}$, $C_i$, $C_{i+1}$, $\mathcal T_i$, 
$\mathcal T_{i+1}$ are as in Proposition~$\ref{prop:conj}$, that $H=C_iL_i$, $K=C_{i+1}L_{i+1}$, 
and that  
$L_{i+1}=\left<a\right>$. 
Then $K=C_{i+1}=N_{G_0}(\la a\ra)$ and
$$
|\mathcal{T}_{i+1}|=|N_{G_0}(\la a\ra):C_{G_0}(a)|\leq \varphi(|a|). 
$$
More\-over,
if $({G_0},L_i,\la a\ra,\ismember)$ satisfies the membership test
condition in $G$, for some algorithm $\ismember$, then so does
$({G_0},H,N_{G_0}(\left<a\right>),\newismember)$ where the algorithm $\newismember$ is given by \mbox{Algorithm}~$\ref{newismember}$.
\end{lemma}

\begin{proof} By the definition of $\mathcal T_{i+1}$ and $L_{i+1}$
\[
a^{\mathcal{T}_{i+1}L_{i+1}}= a^{G_0}\cap L_{i+1}=
a^{G_0}\cap\la a\ra = a^{N_{G_0}(\la a\ra)}.
\] 
However, $L_{i+1}$ centralises $a^{N_{G_0}(\la a\ra)}$ and so $a^{\mathcal{T}_{i+1}}=
a^{N_{G_0}(\la a\ra)}$, which implies that 
$$
C_{i+1}=
C_{G_0}(a)\mathcal{T}_{i+1}=N_{G_0}(\la a\ra).
$$
Moreover, since $L_{i+1}=\la a\ra \leq C_{G_0}(a)\subseteq 
N_{G_0}(\la a\ra)=C_{i+1}$, we obtain that $K=C_{i+1}$.
 
Since $L_{i+1}$ is abelian, $|\mathcal{T}_{i+1}|=|a^{G_0}\cap\la a\ra|$, and since 
$N_{G_0}(\la a\ra)$ acts on the set of $\varphi(|a|)$ 
generators of $\la a\ra$,
with kernel $C_{G_0}(a)$ and with $a^{G_0}\cap L_{i+1}$ as one of the orbits, 
it follows that $|a^{G_0}\cap L_{i+1}|=|a^{N_{G_0}(\la a\ra)}|=|N_{G_0}(\la a\ra):C_{G_0}(a)|$.
The final assertion is part (b) of
Lemma~\ref{lem:conj1}.\qed
\end{proof}

\section{{\sc IsMember} using element orders}\label{sect:orders}

In this section we present a version of {\sc BasicSift} that has proved
useful especially for the first link in a chain such as (\ref{chain}) for several sporadic simple groups
$G$. It requires the relevant subsets 
to be subgroups. We give some applications
that use this version in Section~\ref{sect:appns}.

As in Section~\ref{sect:conj}, we will describe a version of the procedure {\sc IsMember} that can be used in the $\basicsift$ Algorithms~\ref{basicsift} 
and~\ref{basicsiftcr}.
Let $G$ and $G_0$ be finite groups such that
$G_0\leq G$, and 
suppose that $H$ and $K$ are subgroups of $G_0$, with $K<H$. Therefore
condition (\ref{eq:hl}) automatically holds with $L=H$.
An extra requirement is that for all subgroups $M$ such that $K<M\leq H$,  
a reasonable proportion of the elements of
$M$ have orders that do not occur as orders of elements in $K$. 
We define
$$
I=\{n\in\mathbb N\ |\ \mbox{some $M$ with $K<M\leq H$ 
has elements of order $n$ but $K$ does not}\}.
$$
Assume that $I\neq\emptyset$ and let $p_0$ be a number
such that for all $M$ with $K<M\leq H$ the proportion of the elements of $M$ with
orders in $I$ is at least $p_0$. We suppose that $p_0>0$.
As usual we assume 
that random selections in the procedure are made independently 
and uniformly from the relevant subgroups.
Moreover, we emphasise that this is a `black-box algorithm' , and in particular it is not easy to find the order of an element efficiently. To test if an element $g$ has a particular order $n\in I$, we check first that $g^n=1$ which implies that the order of $g$ divides $n$, and then, for each maximal proper divisor $d$ of $n$, we test that $g^d\ne 1$.   
We define $\bar I$ to be the number of integers that are
either equal to or a
maximal proper divisor of an element of $I$. Then for $g\in G_0$ we can
test if the order of $g$ lies in $I$ by examining $\bar I$ powers of $g$.

\begin{algorithm}
\SetKw{Input}{Input:}
\SetKw{Output}{Output:}
\SetKw{Set}{set}
\SetKwFunction{IsMember}{IsMember}
\SetKwFunction{Random}{\sc RandomElement}
\caption{The algorithm {\sc IsMemberOrders}}
\label{ismember-orders}
Algorithm~\ref{ismember-orders}: 
{\sc IsMemberOrders}\\
\tcc{See Proposition~\ref{prop:orders} for notation}
\medskip
\Input{\rm  $(y,e)$ where $y\in G$, and $0< e<1/2$}\;
\Output{\true\ {\rm or} \false}\;
\Set $N=\lceil\log(e^{-1})/\log((1-p_0)^{-1})\rceil$\;
\Set{$n=0$}\;
\Repeat{$n\geq N$}{\Set $h=\Random(\left<K,y\right>)$\;
\If{\rm the order of $h$ is in $I$}{\Return\false}
\Set{$n:=n+1$}}
\Return\true
\end{algorithm}

\begin{proposition}\label{prop:orders} 
Suppose that $G$, $G_0$, 
$H$, $K$, $I$,  $\bar I$, and $p_0$ are as above. Also suppose that, for any $M$ satisfying $K\leq M\leq H$, $\randel(M)$ returns  uniformly distributed, independent random elements of $M$.  
Then $(G_0,H,K,\mbox{\sc IsMember})$ satisfies
the membership test condition in $G$, 
where {\sc IsMember} is Algorithm~$\ref{ismember-orders}$. 
Further, the cost of running $\mbox{\sc IsMemberOrders}(\cdot,e)$ is 
\[
O\big(\log(e^{-1})\cdot p_0^{-1}\,(\xi +  \log(\max I)\cdot\bar{I}\cdot \varrho)\big)
\]
where $\max I$ is the maximum integer in $I$, and
$\varrho$, $\xi$ are upper bounds for the costs of a group
operation in $G$, and making a random selection from
any subgroup of the form $\la K, g\ra$ ($g\in G$), 
respectively.
\end{proposition}
\textbf{Remark:} In Algorithm~\ref{ismember-orders}
we have to make a random selection from a possibly different
group $\langle K,y \rangle$ for every step of the loop. Because the
known algorithms for producing (pseudo-) random elements in groups all
involve an initialisation phase, the constant $\xi$ here could be much bigger
than the constant $\rho$ or even the corresponding constant $\xi$ in
other algorithms of this paper.

\begin{proof} If $y\in K$, then by one of the conditions on the input, no element of 
$\la K,y\ra = K$ has order in $I$, and hence the output is {\sc True}. Now suppose that $y\in H
\setminus K$  so that $K<\la K,y\ra\leq H$. By assumption, the proportion of elements of 
$\la K,y\ra$ with order in $I$ is at least $p_0$. Thus, after $N$ independent random selections 
from $\la K,y\ra$, the probability
that we do not find at least one element with order in $I$ is at most $(1-p_0)^N$. The definition
of $N$ implies that $(1-p_0)^N\leq e$. Thus the membership test
condition is
satisfied.

Now we estimate the cost.
For each random $h\in  \la K, y\ra$, we compute $h^n$ for
each $n$ that is either equal to or a maximal divisor of an element of $I$. We do this by 
first computing $h^2,h^4,\dots,h^{2^m}$, where $2^m\leq{\max I}<2^{m+1}$.
We use these elements to compute $h^n$, for each relevant $n$, with at most $m\bar{I}$ group 
multiplications.
Thus the cost of computing all of the relevant $h^n$ is at most $m\bar{I}\varrho=O(\log({\max I}) \bar{I} 
\varrho)$.
The number of random $h$ to be processed is at most $N$, which, by Lemma~\ref{rem:mc}, is
$O(\log(e^{-1}) \cdot p_0^{-1})$. 
Thus an upper bound for the cost is
$O\big(\log(e^{-1})\cdot p_0^{-1}\,(\xi +   \log({\max I})\cdot\bar{I}\cdot\varrho)\big)$.
\qed
\end{proof}

In most cases when Algorithm~\ref{ismember-orders} is used, we
have that $K$ is maximal in $H$, and so the only possibility for 
$M$ in Proposition~\ref{prop:orders} is $K$ or $H$. 
Also it is often true that $I$ consists entirely of primes, and then
$\bar{I}=|I|+1$.

\begin{corollary}\label{final}
Use the notation of Proposition~$\ref{prop:orders}$ and suppose that $u=|H:K|$. Let
$\basicsift$ be Algorithm~$\ref{basicsift}$ with
Algorithm~$\ref{ismember-orders}$ as $\ismember$. 
Then
the cost of executing $\basicsift(\cdot,\ve)$ with $0<\ve<1/2$ is
\[
O\left( \log(\ve^{-1})\cdot u\left(\xi + \varrho+\frac{\log(\ve^{-1})+\log u}{p_0}\,(\xi' +
\log({\max I})\cdot\bar{I}\cdot\varrho) \right) \right),
\]
where $\xi$ is the cost of selecting a random element of $H$, $\xi'$ 
is an upper bound for the cost
of selecting a random element from a subgroup of the form $\left<K,x\right>$,  
where $x\in H$, and $\varrho$ is the cost of a group operation in $G$.
\end{corollary}
\begin{proof} Using the notation of Theorem~\ref{prop:basic-mc}, since $H=L>K$, we have $p=|K|/|H|$, which is $u^{-1}$. Thus, by Theorem~\ref{prop:basic-mc} and Proposition~\ref{prop:orders}, the cost of this version of 
{\sc BasicSift}$(\cdot,\ve)$ is
\[
O\big( \log(\ve^{-1})\cdot u\big(\xi + \varrho+\log(e^{-1})p_0^{-1}\,(\xi' +
\log({\max I})\cdot\bar{I}\cdot\varrho) \big) \big),
\]
where $e=\ve u^{-1}/2(1-u^{-1})$. Now 
$$
\log(e^{-1}) = \log(\ve^{-1}) +\log(2) + \log(u-1)
= O( \log(\ve^{-1}) +\log u),
$$ 
and the assertion follows. \qed
\end{proof}

\section{The Higman-Sims group HS revisited}
\label{HSrevisited}

In Section~\ref{example} we presented a simple algorithm to write an element
of $\HS$ as a word in a given generating set. This algorithm served as
an example for the theory developed in this paper. We now examine
how the steps of the $\HS$ algorithm in Section~\ref{example} 
fit into the theoretical
framework presented in the subsequent sections. 
We use the notation of Section~\ref{example}.

As in Section~\ref{example}, $G$ is a group isomorphic to $\HS$, and 
we set $G_0=G$.
Let $L_1$ be a maximal subgroup of $G$ isomorphic to 
$U_3(5).2$. Then $L_1$ has a subgroup $Z$ of order~16. We noted in
Section~\ref{example} that the proportion of elements of order
$11$ or $15$ in $\HS$ is $41/165$, while $L_1$ does not contain any such
element. Let $\ismember_1$ be Algorithm~\ref{ismember-orders} with $I=\{11,15\}$ and
$p_0=41/165$. Then, by Proposition~\ref{prop:orders}, 
$(G,G,L_1,\ismember_1)$ satisfies
the membership test condition in $G$. Let $C_1=C_G(a)$ where $a\in Z$ and $|a|=8$
as in Section~\ref{example} and
let $\newismember_1$ 
be Algorithm~\ref{newismember} with $\ismember_1$ as $\ismember$. 
Then, by Lemma~\ref{lem:conj1},
$(G,G,C_1L_1,\newismember_1)$ also satisfies the membership test
condition in $G$, and we use Algorithm~\ref{basicsift} 
to obtain an algorithm $\basicsift_1$
such that $(G,G,C_1L_1,\basicsift_1)$ satisfies the basic sift
condition in  $G$. 

In the next step we recall that $L_2=5^{1+2}:(8:2)$.
We noted that $L_2=N_G(Z(5^{1+2}))$, and so it is
easy to design a deterministic algorithm $\ismember_2$ such that
the $4$-tuple $(G,L_1,L_2,\ismember_2)$ satisfies the membership test
condition in $G$ (just check whether a generator for $Z(5^{1+2})$ is mapped into $Z(5^{1+2})$). 
We set $C_2 = C_G(a) \mathcal{T}_2$ as in Section \ref{sect:secondsubset}. 

Using Algorithm~\ref{newismember}, we find an algorithm 
$\newismember_2$, using $\ismember_2$ as $\ismember$,
such that $(G,C_1L_1,C_2L_2,\newismember_2)$ also satisfies
the membership test condition in $G$, and we use Algorithm~\ref{basicsift} to build an
algorithm $\basicsift_2$ so that 
$(G,C_1L_1,C_2L_2,\basicsift_2)$ satisfies the
basic sift condition in~$G$. 

As $L_3$ is a cyclic group of order~$8$ and $C_3 = C_G(a) \mathcal{T}_3$ as
in Section \ref{sect:thirdsubset}, 
it is easy to check membership
in $L_3$, and following the procedure explained above, it is easy to
obtain an algorithm $\basicsift_3$ such that
$(G,C_2L_2,C_3L_3,\basicsift_3)$ satisfies the basic sift condition 
in~$G$.
In Section~\ref{example} we set $C_4=C_G(a)$, and, using this fact, we
can easily test membership in $C_4$. Thus 
the $4$-tuple $(G,C_3L_3,C_4,\basicsift_4)$ can be constructed.

Finally, it is possible to list all 16 elements of $C_4$ and,
via an exhaustive search, to construct an algorithm $\basicsift_5$
such that $(G,C_4,\{1\},\basicsift_5)$ satisfies the basic sift condition
in~$G$.

Algorithm~\ref{sift} can be used with 
$(G,G,C_1L_1,\basicsift_1)$, $(G,C_1L_1,C_2L_2,\basicsift_2)$,
$(G,C_2L_2,C_3L_3,\basicsift_3)$, $(G,C_3L_3,C_4,\basicsift_4)$,
and $(G,C_4,\{1\},\basicsift_5)$ to sift an element through the chain
$$
G\supset C_1L_1\supset C_2L_2\supset C_3L_3\supset C_4\supset\{1\}.
$$

 \section{Application of the results to sporadic simple groups}
 \label{sect:appns}

An important part of the research presented here is to find explicitly a
suitable subset chain~\eqref{chain} and a $\basicsift$ algorithm
for each step in this chain for many sporadic simple groups. 

Note that all example chains in this section provide pure black-box
algorithms. No particular prior knowledge about the representations
of the groups is used during the sifting. Of course, to construct
the chains we made heavy use of lots of available information and especially
of nice representations.

In the implementations, all occurring group elements are expressed as
straight line programs in terms of standard generators in the sense of
\cite{Wilson96} and \cite{wwwatlas}.

One could improve the performance by using specially crafted
$\ismember$ tests relying on 
specific information about the given representation.
Also, other methods will be better for certain representations.

In this section we assume that $G=G_0$ is one of the sporadic simple groups.
For each group $G$ a subset $S_i$ in the
chain~\eqref{chain} will be a product
$S_i=C_G(a)\mathcal T_iL_i$
with suitable $a$, $\mathcal T_i$, and 
$L_i$. We also set $C_i=C_G(a)\mathcal T_i$ and the 
sequence $C_1,\ldots,C_{k-1}$ will be referred to as a {\it $C$-sequence}.
The ingredients $a$, $L_i$, $\mathcal T_i$ are  in the tables
below. In order to present the subset chains 
in the most compact form, we use the following notation.

{\bf The $a$-column.} 
If the function 
{\newismember} is used to sift through this step of the subset chain,
then this column specifies the conjugacy class of $a$ used by
{\newismember}. The conjugacy class is given using the 
{\sc Atlas} notation; see \cite{Atlas}. 
We can assume without loss of generality that $a$ is contained
in all subgroups $L_i$ where we need the hypothesis $a^G \cap L_i \neq
\emptyset$. If the function 
{\newismember} is not used in this step of the chain then a dash is displayed 
in the appropriate cell.

{\bf The $C_G(a)$-column.} This column contains information
about the centralisers occurring in
the $C$-sequence $C_1,\dots,C_{k-1}$. Note that the $C_i$ 
satisfy the conditions in (\ref{Cchain}).

{\bf The $|\mathcal T_i|$-column.}
Here we only specify the number of elements in $\mathcal{T}_i$. 
In each of the examples, we set $\mathcal T_0=\{1\}$ and, for $i\geq 0$,
the subset $\mathcal
T_{i+1}$ is constructed using the procedure
at the
beginning of Section~\ref{sect:conj}.

{\bf The $L_i$-column.} 
In each table we list the subgroups
$L_1,\ldots,L_{k-1}$ that are used to construct the subgroup chain~\eqref{Lchain}; this chain will be referred to as  the $L$-chain. 
Each
such subgroup is specified as precisely as necessary to define
the descending subset chain. For example, in $\HS$, the group $L_1$ is specified
as $U_3(5).2$ ({\sc Atlas} notation, see \cite{Atlas}), which means that 
any subgroup of $G$ that is isomorphic to $U_3(5).2$ can
play the r\^ole of $L_1$. Similarly, one may
take $L_2$ to be any subgroup  of $L_1$ that is the
semidirect product of an extraspecial group of order $125$ and a 2-group, as
explained in the corresponding cell of the table.

{\bf The $p$-column.} In this column we display the sifting parameter
$p(C_{i-1}L_{i-1},C_{i}L_{i},L_{i-1})$ (see Definition \ref{def:hl} and Proposition \ref{prop:conj}).

{\bf The $\basicsift$-column (BS).} 
We describe the
$\basicsift$ algorithm that is used in a particular step of the subset
chain. The letter R stands for {\basicsiftrandom} (see Algorithm
\ref{basicsift}) and the letter C stands for {\basicsiftcosetreps} (see
Algorithm \ref{basicsiftcr}). Note that in some cases Algorithm
\ref{basicsiftcr} is also used to try a certain set of group elements,
such as the set $\mathcal T_i$ or its inverses.

{\bf The $\ismember$-column.} In this column we describe, how we test
membership in the subgroup $L_i$. If an $a$ is specified in the $a$-column,
then we first design an algorithm $\ismember$ for the pair
$(L_{i-1},L_i)$ using the parameters  in the same cell of the
table.  Then we use Algorithm~\ref{newismember} to obtain a new
algorithm $\ismember$ for the pair $(C_{i-1}L_{i-1},C_iL_i)$, and
finally,  Algorithm~\ref{basicsift} yields a $4$-tuple
$(G,C_{i-1}L_{i-1},C_iL_i,\basicsift_i)$ satisfying the basic sift condition
in~$G$. 

The membership test $\ismember$ for the pair $(L_{i-1},L_i)$ is described
using the following notation.

(a) If a set $I$ of element orders is specified, Algorithm~\ref{ismember-orders}
is used for the $\ismember$ test for $L_i$. In this case we also specify the
probability $p_0$ to find an element of such an order in $L_{i-1}$.

(b) If, in
the $\basicsift$-column of the table, an $L_i$ is specified to
be the centraliser or the normaliser of an element or a subgroup, then, using this
fact,  we build a
deterministic algorithm to determine membership of $L_i$.

(c) Finally, the symbol 
$1$ in that column indicates that we use an exhaustive search to 
test equality in the subgroup $L_i$. This method will be used in 
the special case when $L_i=1$. 

Note that the symbol ``1'' may stand either for the trivial subgroup 
or for the identity element, but its meaning is always clear from the context.

\begin{table}
$$
\begin{array}{cccccccc}
\hline
M_{11}& a       & C_G(a)         &|\mathcal{T}_i|& L_i      & p  &\mbox{\sc BS}
     &\ismember \\
\hline
  1  & a\in2A  & 2.S_4          &  2  & 2.S_4    &13/165 & R 
     & C_G(a) \\
  2  & a\in2A  & 2.S_4          &  3  & 2^2      & 1/6 & C
     & C_{L_1}(b) \mbox{ with } b^2=1 \\
3 & a\in 2A & 2.S_4 & 1 & 1 & 1/3 & C & 1\\
\hline
  3  &    -    & -              &  -  & 2.S_4    & - & -
     & - \\
  4  &    -    & -              &  1  & 8        & 1/6 & C
     & C_{L_3}(8A) \\
  5  &    -    & -              &  1  & 1        & 1/8 & C
     & 1 \\
\hline
\end{array}
$$
\caption{A chain for $M_{11}$ using $2.S_4$}
\label{M11_1}
\end{table}

The first example is in Table \ref{M11_1}, which describes a subset
chain for the sporadic simple Mathieu group $M_{11}$. 
In Table \ref{M11_2} we present another subset chain for $M_{11}$ to
demonstrate a new idea, namely that information gained during an {\sc
IsMember} test can be used further.
Table \ref{M12} contains a subset chain for the sporadic simple
Mathieu group $M_{12}$.
In Table \ref{M22} we describe a subset chain for the sporadic simple
Mathieu group $M_{22}$.
Table \ref{J2_1} presents a subset chain for the sporadic simple
Janko group $J_2$, that uses only deterministic membership tests.
In contrast, Table \ref{J2_2} shows another chain for $J_2$ with
membership tests using element orders.

We conclude this section with a larger example, in which we demonstrate
yet another idea, namely that there may be ``branches'' in chains,
leading to different behaviour of the algorithm under certain circumstances,
that may occur during the calculation. See Table \ref{Ly} for details 
and Note~(i) to Table~\ref{Ly} for an explanation.

We have implemented the generalised sifting algorithms using the 
subset chains described in the tables below for some of the sporadic 
simple groups. The implementations were written in the~{\sf GAP~4} computational
algebra system~\cite{GAP} and will be made available separately 
in the future. Information
on the performance of our implementations can be found in Table~\ref{perf}
and in the notes to that table.

In practical implementations the sifting is carried out in 
several stages. In the first stage we sift our element into a 
smaller subgroup (usually a centraliser of an element), and 
then we start a new sifting procedure in that subgroup. We 
repeat this until we reach the trivial subgroup containing 
only the identity element. In our tables we indicate the 
boundary between different stages by a horizontal line. For 
instance in Table 1, we first sift our element into the 
subgroup 2.S4, and then carry out a new sifting procedure in $2.S_4$.

\begin{table}
$$
\begin{array}{cccccccc}
\hline
M_{11}& a       & C_G(a)         &|\mathcal{T}_i|& L_i      & p  &\mbox{\sc BS}
     &\ismember \\
\hline
  1  & a\in 11A  & \left<a\right>          &  1  & L_2(11)    & 1/12 & C 
     & \mbox{see notes below} \\
  2  & a\in 11A  & \left<a\right>          &  1  & |N_G(\left<a\right>)|=55   & 1/12 & C
& N_G(\left<a\right>) \\
3 & a\in 11A & \left<a\right> & 1 & 1 & 1/5 & C & 1 \\
\hline
  3  &    -    & -              &  -  & \left<a\right>    & - & -
     & - \\
  4  &    -    & -              &  1  & 1    & 1/11 & C
     & 1 \\
\hline
\end{array}
$$
\caption{A second chain for $M_{11}$ using $L_2(11)$}
\label{M11_2}
\end{table}

\subsection*{Notes to Table \ref{M11_2}}

Let $a$ be as in the table and select $x\in G$. We want to write
the element $x$ as a word in a given nice generating set.
Choose an element $a'\in 11A\cap L_1$ such that
$[a,a']\neq 1$ and let $z\in L_1$ with $(a')^z=a$. Then $L_1$ has 12 Sylow 11-subgroups, namely
$\left<a\right>$ and $\left<(a')^{a^i}\right>$ for $i=0,\ldots,10$. For $y_1\in G$, $a^{xy_1}\in L_1$ if and only if $\left<a^{xy_1}\right>$ coincides with one
of the Sylow 11-subgroups of $L_1$. Further, 
such a Sylow subgroup is self-centralising in $G$.
Thus the membership test $a^{xy_1}\in L_1$
is carried out by checking whether $[a^{xy_1},a]=1$
or $[a^{xy_1},(a')^{a^i}]=1$ for some $i\in\{0,\ldots,10\}$. 

The second step of the sifting can be made more efficient as 
follows. Assume that $a^{xy_1}\in L_1$. If $[a^{xy_1},a]=1$ then 
$a^{xy_1}\in L_2$, and we can proceed to the third step of the sifting
procedure. If $[a^{xy_1},(a')^{a^i}]=1$ then, for $y_2=a^{11-i}z$ we have 
that $a^{xy_1y_2}\in L_2$. Thus, storing some information about the
membership test in the first step, we can immediately select the sifting
element $y_2$ in the second step.

\begin{table}
$$
\begin{array}{cccccccc}
M_{12}& a       & C_G(a)         &|\mathcal{T}_i|& L_i      & p  &\mbox{\sc BS}
     &\ismember \\
\hline
 1   & a\in2A   & 2\times S_5    &  1  & M_8.S_4  &1/33& R
     & C_G(2B) \\
 2   & a\in2A   & 2\times S_5    &  1  & |C_{L_1}(x)|=32 & 1/3 & C
     & C_{L_1}(x) \mbox{ with } x^4=1 \\
 3   & a\in2A   & 2\times S_5    &  2  & |C_{L_1}(y)|=8  & 1/2 & C
     & C_{L_1}(y) \mbox{ with } y^4=1 \\
4 & a\in 2A & 2\times S_5 & 1 & 1 & 1/2 & C & C_G(2A) \\
\hline
 4   &    -     &      -         &  -  & 2\times S_5         & - & -
     & - \\
 5   &    -     &      -         &  1  & |N_{L_5}(z)|=40 & 1/6 & C
     & N_{L_5}(z) \mbox{ with } z^5=1 \\
 6   &    -     &      -         &  1  & |C_{L_5}(z)|=10 & 1/4 & C
     & C_{L_5}(z) \mbox{ with } z^5=1 \\
 7   &    -     &      -         &  1  &  1              & 1/10 & C
     & 1 \\
\hline
\end{array}
$$
\caption{A chain for $M_{12}$}
\label{M12}
\end{table}

\begin{table}
$$
\begin{array}{cccccccc}
\hline
M_{22}& a       & C_G(a)         &|\mathcal{T}_i|& L_i      & p  &\mbox{\sc BS}
     &\ismember \\
\hline
  1   & a \in2A & 2^4:S_4        &  1  & L_3(4)   &3/11& R
     & I=\{6,8,11\}, p_0 = 103/364  \\
  2   & a \in2A & 2^4:S_4        &  2  & 2^4:A_5  &5/21& R
     & I=\{7\}, p_0 = 2/7  \\
3 & a\in 2A & 2^4:S_4 & 1 & 1 & 1/60 & R & C_G(a)\mbox{ see Note below}\\
\hline
  3   &    -    &   -            &  -  & 2^4:S_4  &-& -
     & - \\
  4   &    -    &   -            &  1  & 2^4      &1/24& C
     & C_G(x) \mbox{ with } z^2=1  \\
  5   &    -    &   -            &  1  &  1       &1/16& C
     & 1  \\
\hline
\end{array}
$$
\caption{A chain for $M_{22}$}
\label{M22}
\end{table}

\subsection*{Notes to Table~\ref{M22}}

Here the elements from $\mathcal{T}_2 = \{ 1,t_1 \}$ are tried
together with elements from the group $L_2$ to reach the centraliser
of $a$. The probability $1/60$ is the minimum of the probability for
the two cases $C_G(a) \cdot \{1\} \cdot L_2$ and $C_G(a) \cdot \{t_1\}
\cdot L_2$.

\begin{table}
$$
\begin{array}{cccccccc}
\hline
 J_2 & a       & C_G(a)         &|\mathcal{T}_i|& L_i      & p  &\mbox{\sc BS}
     &\ismember \\
\hline
   1 & a\in8A  & \left<a\right> &  2  & 3.A_6.2_2&1/140& R
     & N_G(3A)=N_G(\mbox{Soc}(L_1))  \\
   2 & a\in8A  & \left<a\right> &  4  & 3^{1+2}:8& 1/5& C
     & N_G(3^{1+2})=N_G(\mbox{Syl}_3(L_1))  \\
   3 & a\in8A  & \left<a\right> &  4  & 8        &1/27& C
     & C_G(a) = \left< a \right>  \\
4 & a\in 8A & \left<a\right> & 1 & 1 & 1/4 & C & C_G(a) \\
\hline
   4 &   -     & -              &  -  & \left<a\right>        &- & -
     &-  \\
   5 &   -     & -              &  1  & 1        &1/8 & C
     & 1        \\
\hline
\end{array}
$$
\caption{A chain for $J_2$ with deterministic membership tests}
\label{J2_1}
\end{table}

\begin{table}
$$
\begin{array}{ccccccccc}
\hline
 J_2 & a       & C_G(a)         &|\mathcal{T}_i|& L_i      & p  &\mbox{\sc BS}
     &\ismember& \mbox{Note} \\
\hline
   1 & a\in2A  & 2_-^{1+4}:A_5  &  1  & 3.A_6.2_2 &1/6 &  R     
     & N_G(3A) & \mbox{(i)} \\
   2 & a\in2A  & 2_-^{1+4}:A_5  &  1  & 3\times A_5    &1/3 & C 
     & I=\{4,12\}, p_0 = 1/4 & \mbox{(ii)} \\
   3 & a\in2A  & 2_-^{1+4}:A_5  &  1  & A_4            &1/5 & C
     & I=\{5\}, p_0 = 2/5 & \\
   4 & a\in2A  & 2_-^{1+4}:A_5  &  1  & 4              &1/3 & C
     & C_{L_3}(a) & \\
\hline
4 & - & - & - & 2^{1+4}_-:A_5 & - & - & - \\
   5 &    -    &         -      &  1  & |L_5| = 192    &1/10& C
     & C_{C_1}(2C) & \\
   6 &    -    &         -      &  1  & |L_6| = 32     &1/6 & C
     & N_{C_1}(4A) & \\
   7 &    -    &         -      &  1  & |L_7| = 16     &1/2 & C
     & C_{C_1}(4A) & \\
   8 &    -    &         -      &  1  & 1              &1/16& C
     &    1    &  \\
\hline 
\end{array}
$$
\caption{Another chain for $J_2$}
\label{J2_2}
\end{table}

\subsection*{Notes to Table \ref{J2_2}}
\begin{itemize}
\item[(i)] $a^G \cap (3.A_6.2) \leq 3.A_6$, so we get an index $2$ for free.
\item[(ii)] The $3$ of $3 \times A_5$ is in $C_G(a)$, 
and hence $C_2L_2=2\times A_5$.
\end{itemize}

\begin{table}
$$
\begin{array}{cccccccc}
\hline
 \HS & a & C_G(a) &|\mathcal{T}_i|& L_i  & p  &\mbox{\sc BS}&\ismember\\
\hline
   1 & a\in 8B & 2 \times 8     &   1 & U_3(5).2 &1/88&   R    & 
       I=\{11,15\}, p_0 = 41/165 \\
   2 & a\in 8B & 2 \times 8     &   2 & 5^{1+2}:(8:2) & 1/63& R&
       N_G(Z(5^{1+2})) \\
   3 & a\in 8B & 2 \times 8     &   4 & \left<a\right>& 1/125&R&
       C_{L_1}(a) \\
   4 & a\in 8B & 2 \times 8     &   1 & 1             & 1/4  &C&
       1      \\
\hline
4 & - & - & - & 2\times 8 & - & - & - \\
   5 &    -    &       -        &   1 & 1             & 1/16 &C&
       1      \\
\hline 
\end{array}
$$
\caption{The chain for $\HS$ from sections \ref{example} and
\ref{HSrevisited}}
\label{HS_1}
\end{table}

In Sections \ref{example} and \ref{HSrevisited} we already described the
subgroup chain for the sporadic simple Higman-Sims group $\HS$ presented
in Table \ref{HS_1}. We
found this chain very useful to illustrate the ideas used in this paper.
However, it turns out that one can design a much more efficient chain
for $\HS$ whose details are presented in Table \ref{HS_2}.

\begin{table}
$$
\begin{array}{cccccccc}
\hline
 \HS & a & C_G(a) &|\mathcal{T}_i|& L_i  & p  &\mbox{\sc BS}&\ismember\\
\hline
  1 & a \in 2A & 4.2^4:S_5 & 1 & M_{22} & 1/5  & R &
      I=\{ 10,12,15,20 \}, p_0=7/20 \\
  2 & a \in 2A & 4.2^4:S_5 & 1 & L_3(4) & 3/11 & R &
      I=\{ 6,8,11 \}, p_0 = 103/264 \\
  3 & a \in 2A & 4.2^4:S_5 & 1 & A_6 & 1/7 & R &
      I=\{ 7 \}, p_0 = 2/7 \\
  4 & a \in 2A & 4.2^4:S_5 & 1 & A_5 & 1/3 & C &
      I=\{ 4 \}, p_0 = 1/4 \\
  5 & a \in 2A & 4.2^4:S_5 & 1 & A_4 & 1/5 & C &
      \mbox{(i)} \\
  6 & a \in 2A & 4.2^4:S_5 & 1 & 2^2 & 1/3 & C &
      C := C_G(a) \mbox{ (ii)} \\
\hline
6 & - & - & - & 4.2^4:S_5 & - & - & - \\
  7 & -        & -         & 1 & 4.2^4:A_5 & 1/2 & C &
      C_C(4B) \\
  8 & -        & -         & 1 & 4.2^4.2^2 & 1/15 & C &
      N_C(x^2) \mbox{ for some } x \mbox{ with } x^8=1 \\
  9 & -        & -         & 1 & 4.2^4.2   & 1/2  & C &
      C_C(x^2) \\
 10 & -        & -         & 1 & 8 \times 2 & 1/8 & C &
      C_C(x) \\
 11 & -        & -         & 1 & 1          & 1/16 & C &
      1      \\
\hline
\end{array}
$$
\caption{More efficient chain for $\HS$}
\label{HS_2}
\end{table}

\eject

\subsection*{Notes to Table \ref{HS_2}} 
\begin{itemize}
\item[(i)] The $2^2$ in $A_4$ is equal to $C_{A_4}(a)$, therefore
we can test membership of $a^g$ in $A_4$ efficiently.
\item[(ii)] Here we reach $C_G(a)$, since $2^2 \leq C_G(a)$.
\end{itemize}

\begin{table}
$$
\begin{array}{ccccccccc}
\hline
 \Ly & a  & C_G(a) &|\mathcal{T}_i|& L_i & p  &\mbox{\sc BS}&\ismember &
\mbox{Note} \\
\hline
  1 & a\in 3A & 3.McL    &   3 & 3.McL  & 15401/9606125 & R &
      C_G(a) & \mbox{(i)} \\
  2 & a\in 3A & 3.McL    &   1 & 2.A_8  & 1/275 & R &
      C_{L_1}(2A) & \mbox{(ii)} \\
  3 & a\in 3A & 3.McL    &   3 & 3 \times (2.A_5) & 11/56 & R &
      C_{L_2}(3A) & \mbox{(iii)} \\
  4 & a\in 3A & 3.McL    &   3 & 3 \times (2.S_3) & 1/10 & C &
      N_{L_3}(3B) & \\
  5 & a\in 3A & 3.McL    &   4 & 3 \times 2 \times 3 & 1/2 & C &
      \mbox{set} & \mbox{(iv)} \\
\hline
5 & - & - & - & 3.{McL} & - & - & -\\
  6 & a'\in 3C & 3^{2+4}.(2.A_5) & 1 & 2.A_8 & 1/275 & R &
      C_{3.McL}(2A) & \mbox{(v)} \\
  7 & a'\in 3C & 3^{2+4}.(2.A_5) & 3 & 3 \times (2.A_5) & 11/56 & R &
      C_{L_6}(3A) & \\
  8 & a'\in 3C & 3^{2+4}.(2.A_5) & 3 & 3 \times (2.S_3) & 1/10 & C &
      N_{L_7}(3B) & \\
  9 & a'\in 3C & 3^{2+4}.(2.A_5) & 4 & 3 \times 2 \times 3 & 1/2 & C &
      \mbox{set} & \mbox{(vi)} \\
\hline
9 & - & - & - & 3^{2+4}.(2.A_5) & - & - & - \\
 10 & -        & -               & 1 & |L_{10}|=1080 & 1/81 & C &
      C_{C'}(z) & \mbox{(vii)} \\
 11 & -        & -               & 1 & |L_{11}|=90  & 1/12 & C &
      C_{C'}(z') & \mbox{(viii)} \\
 12 & -        & -               & 1 & |L_{12}|=9   & 1/10 & C &
      \mbox{Syl}_3(L_{11}) & \mbox{(ix)} \\
 13 & -        & -               & 1 & 1            & 1/9 & C &
      1 & \\
\hline
\end{array}
$$
\caption{A chain for $\Ly$}
\label{Ly}
\end{table}

\subsection*{Notes to Table \ref{Ly}}
\begin{itemize}
\item[(i)] $a^G \cap L_1 = \{ a,a^{-1} \} \cup a^{xL_1}$ for some $x \in G$.
We store an element $y \in G$ with $a^y = a^{-1}$ and handle the cases
$a^g = a$ and $a^g = a^{-1}$ separately, which allows us to jump directly
to step $6$ in these cases. Otherwise, we can work with a single
conjugacy class $a^{xL_1}$ in $L_1$. Of course, most of the time this
latter case will occur, as $a^{xL_1}$ has $30800$ elements.
\item[(ii)] The centraliser of a $2A$ element in $3.McL$ is 
$3 \times (2.A_8)$. However, we already have avoided the $3$ in the center
by the special cases in step $1$. Note that we have reduced the number 
$|\mathcal{T}_2|$ to $1$, because $\mathcal{T}_2 = \{x\}$, again by the
special cases in step $1$.
\item[(iii)] $L_3$ is the centraliser in $L_2$ of an element of order $3$.
\item[(iv)] In this step we store the complete set of $4$ possible results
for $a^g$ together with elements of $G$ to conjugate them back to $a$.
So we can reach $C_G(a)$ after this step with no additional costs.
\item[(v)] $a'$ is from $3C$ in $3.McL$. By $3^{2+4}$ in $L_6$ we mean
a $3$-group with an elementary-abelian center of order $9$ with an 
elementary-abelian group of order $81$ as factor group. As in (ii) is the
centraliser of a $2A$ element in $3.McL$ is $3 \times (2.A_8)$.
However, since $a'$ lies in $3C$ of $3.McL$, we automatically reach
$2.A_8$.
\item[(vi)] Note (iv) applies analogously.
\item[(vii)] $z$ is an involution in $C' := C_{3.McL}(a') = 3^{2+4}.(2.A_5)$.
\item[(viii)] $z'$ is an element of order $15$ in $C'$.
\item[(ix)] The Sylow-$3$-subgroup is normal, therefore just looking for
element orders tests membership.
\end{itemize}

\subsection*{Notes to Table~\ref{perf}}
The algorithms presented in this paper
were implemented for the sporadic simple groups above. We used matrix 
representations of these groups and Table~\ref{perf} contains some average
running times in seconds. For each representation, we sifted $1000$
pseudo random elements and the running times are for those $1000$
calls to {\sc Sift} on a machine with a Pentium IV processor running
at $2.53$ GHz with $512$ MB of main memory. The third column contains
the average number of multiplications necessary for one call to {\sc
Sift}, including the generation of pseudo random elements. Note that
the initialization phase of the pseudo random generator (using product
replacement) involves $100$
multiplications for every newly generated group object. In all cases
the bound for the error probability was $1/100$.

\begin{table}\label{perf}
\[\begin{array}{|c|c|c|}
\hline
\mbox{Group} 
  & \mbox{Time for }1000\mbox{ calls in seconds}
  & \mbox{Av.~no.~of mults. per call} \\
\hline
\hline
M_{11}\leq \gl{10}{2} & 0.8 & 116 \mbox{ (chain in Table \ref{M11_1})} \\
\hline
M_{11}\leq \gl{10}{2} & 0.7 & 187 \mbox{ (chain in Table \ref{M11_2})} \\
\hline
M_{11}\leq \gl{45}{3} & 26.1 & 116 \mbox{ (chain in Table \ref{M11_1})} \\
\hline
M_{11}\leq \gl{45}{3} & 34.1 & 185 \mbox{ (chain in Table \ref{M11_2})} \\
\hline
M_{12}\leq \gl{10}{2} & 2.2 & 243 \\
\hline
M_{12}\leq \gl{44}{2} & 11.6 & 245 \\
\hline
M_{12}\leq \gl{16}{11} & 10.9 & 238 \\
\hline
M_{22}\leq \gl{10}{2} & 25.2 & 1670 \\
\hline
M_{22}\leq \gl{21}{3} & 98.8 & 1685 \\
\hline
J_{2}\leq \gl{36}{2} & 34.6 & 908 \mbox{ (chain in Table \ref{J2_1})} \\
\hline
J_{2}\leq \gl{36}{2} & 46.2 & 847 \mbox{ (chain in Table \ref{J2_2})} \\
\hline
J_{2}\leq \gl{14}{5} & 32.7 & 923 \mbox{ (chain in Table \ref{J2_1})} \\
\hline
J_{2}\leq \gl{14}{5} & 33.2 & 846 \mbox{ (chain in Table \ref{J2_2})} \\
\hline
\HS\leq \gl{20}{2}  & 344.8 & 13923 \mbox{ (chain in Table \ref{HS_1}) } \\
\hline
\HS\leq \gl{20}{2}  & 77.7 & 2783 \mbox{ (chain in Table \ref{HS_2})}\\
\hline
\HS\leq \gl{49}{3}  & 699.2 & 2807 \mbox{ (chain in Table \ref{HS_2})}\\
\hline
\Ly\leq \gl{111}{5} & 19835.8 & 7416 \\
\hline
\end{array}
\]
\caption{Timings and number of multiplications for various chains}
\end{table}

\section*{Acknowledgment}
The research presented in this paper 
forms part of the first author's PhD project, 
supported by an Australian Postgraduate award, 
and was also funded by the Australian Research Council Discovery Grant 
DP0557587. Much of the work leading to this article was carried out
while the fourth author was employed as a Research Associate in 
the Department of Mathematics and Statistics of The University of Western
Australia; he was also supported by the Hungarian Scientific 
Research Fund (OTKA)
grants F049040 and T042706.

We wish to express our thanks to Eamonn O'Brien for helpful comments on an
earlier version of this paper, and also to an anonymous referee 
for a number of perceptive 
observations and queries on our submitted draft, which, in each instance,
led to an improvement in the article. 

\def\cprime{$'$} \def\cprime{$'$} \def\cprime{$'$}


\begin{thebibliography}{10}

\bibitem{BealsLeedhamGreenetal02}
Robert M. Beals, Charles~R. Leedham-Green, Alice~C. Niemeyer, Cheryl~E. Praeger,
  and {\'A}kos Seress.
\newblock Permutations with restricted cycle structure and an algorithmic
  application.
\newblock {\em Combin. Probab. Comput.}, 11(5):447--464, 2002.

\bibitem{BealsLeedhamGreenetal03}
Robert~M. Beals, Charles~R.\ Leedham-Green, Alice~C.\ Niemeyer, Cheryl~E.\
  Praeger, and \'Akos Seress.
\newblock A black-box group algorithm for recognizing finite symmetric and
  alternating groups {I}.
\newblock {\em Trans. Amer. Math. Soc.}, 355(5):2097--2113, 2003.

\bibitem{Magma}
Wieb Bosma, John Cannon, and Catherine Playoust.
\newblock The {M}agma algebra system {I}: {T}he user language.
\newblock {\em J. Symbolic Comput.} 24(3--4):235--265, 1997.



\bibitem{BratusPak00}
Sergey Bratus and Igor Pak.
\newblock Fast constructive recognition of a black box group isomorphic to
  ${S}\sb n$ or ${A}\sb n$ using {G}oldbach's conjecture.
\newblock {\em J. Symbolic Comput.}, 29(1):33--57, 2000.

\bibitem{CellerLeedhamGreen98}
F.~Celler and C.~R. Leedham-Green.
\newblock A constructive recognition algorithm for the special linear group.
\newblock In {\em The atlas of finite groups: ten years on (Birmingham, 1995)},
  pages 11--26. Cambridge University Press, Cambridge, 1998.

\bibitem{CoopermanFinkelsteinetal97a}
Gene Cooperman, Larry Finkelstein, and Steve Linton.
\newblock Constructive recognition of a black box group isomorphic to ${\rm
  {G}{L}}(n,2)$.
\newblock In {\em Groups and computation II (New Brunswick, NJ, 1995)}, pages
  85--100. Amer. Math. Soc., Providence, RI, 1997.




\bibitem{Atlas}
J.~H. Conway, R.~T. Curtis, S.~P. Norton, R.~A. Parker, and R.~A. Wilson.
\newblock {\em Atlas of finite groups}.
\newblock Oxford University Press, Eynsham, 1985.

\bibitem{GAP}
The GAP~Group, \emph{ {GAP -- Groups, Algorithms, and Programming, Version
{\rm 4.4}}}
(Aachen, St~Andrews, 2004), {\tt
http://www.gap-system.org}.

\bibitem{hlo'brw}
P. E. Holmes, S. A. Linton, E. A. O'Brien, A. J. E. Ryba and R. A. Wilson.
\newblock Constructive membership testing in black-box groups. 
\newblock Unpublished manuscript, 2004.

\bibitem{KantorSeress01}
William~M. Kantor and {\'A}kos Seress.
\newblock Black box classical groups.
\newblock {\em Mem. Amer. Math. Soc.}, 149(708):viii+168, 2001.

\bibitem{LeedhamGreen00}
Charles~R.\ Leedham-Green.
\newblock The computational matrix group project.
\newblock In William~M. Kantor and {\'A}kos Seress, editors, {\em Groups and
  {C}omputation {I}{I}{I}}, pages 85--101. {O}{S}{U}
  {M}athematical {R}esearch {I}nstitute {P}ublications, Walter de {G}ruyter,
  2000.

\bibitem{LeedhamGreenNiemeyeretal03}
Charles~R. Leedham-Green, Alice~C. Niemeyer, E.A. O'Brien, and Cheryl~E.
  Praeger.
\newblock Recognising matrix groups over finite fields.
\newblock In V.~Weispfenning J.~Grabmeier, E.~Kaltofen, editors, {\em Computer
  Algebra Handbook, Foundations, Applications, Systems}, pages 459--460,
  Springer-Verlag, Berlin, New York, 2003.

\bibitem{rob}
Derek J. S. Robinson. {\em A course in the theory of groups.}
\newblock Springer-Verlag, 1982.

\bibitem{Seress}
\'Akos Seress. 
\newblock {\em Permutation group algorithms.}
\newblock Cambridge University Press, 2003.

\bibitem{Sims70}
Charles~C.\ Sims.
\newblock Computational methods in the study of permutation groups.
\newblock In {\em Computational problems in abstract algebra}, pages 169--183,
  Oxford, 1970. (Oxford, 1967), Pergamon Press.

\bibitem{Sims97}
Charles~C. Sims.
\newblock Computing with subgroups of automorphism groups of finite groups.
\newblock In {\em Proceedings of the 1997 International Symposium on Symbolic
  and Algebraic Computation (Kihei, HI)}, pages 400--403 (electronic), New
  York, 1997. ACM.

\bibitem{Wilson96} R.~A.~Wilson. \newblock
Standard generators for sporadic simple groups.
\newblock {\em J. Algebra }{\bf 184} (1996), no.~2, 505--515.

\bibitem{wwwatlas} Robert Wilson et al. Atlas of finite group
representations, available on the Internet at \\ 
{\tt http://for.mat.bham.ac.uk/atlas/v2.0/}

\end{thebibliography}
\end{document}